\documentclass[12pt]{amsart}

\usepackage{amssymb}
\usepackage{amscd}

\newcommand{\bt}{\begin{theorem}}
\newcommand{\et}{\end{theorem}}
\newcommand{\bp}{\begin{proposition}}
\newcommand{\ep}{\end{proposition}}
\newcommand{\bq}{\begin{question}}
\newcommand{\eq}{\end{question}}
\newcommand{\bl}{\begin{lemma}}
\newcommand{\el}{\end{lemma}}

\newcommand{\br}{\begin{result}}
\newcommand{\er}{\end{result}}
\newcommand{\be}{\begin{equation}}
\newcommand{\ee}{\end{equation}}
\newcommand{\bc}{\begin{corollary}}
\newcommand{\ec}{\end{corollary}}
\newcommand{\bex}{\begin{example}}
\newcommand{\eex}{\end{example}}

\newtheorem{theorem}{Theorem}[section]
\newtheorem{corollary}[theorem]{Corollary}
\newtheorem{lemma}[theorem]{Lemma}
\newtheorem{proposition}[theorem]{Proposition}
\newtheorem{result}[theorem]{Result}
\newtheorem{example}[theorem]{Example}
\newtheorem{question}[theorem]{Question}

\numberwithin{equation}{section}

\newcommand{\N}{\mathbb{N}}

\newcommand{\Z}{\mathbb{Z}}

\newcommand{\cB}{\mathcal{B}}

\newcommand{\cH}{\mathcal{H}}

\newcommand{\cK}{\mathcal{K}}
\newcommand{\cL}{\mathcal{L}}

\newcommand{\cO}{\mathcal{O}}
\newcommand{\cR}{\mathcal{R}}

\newcommand{\cSb}{\mathcal{S}\!\mathit{ub}}

\newcommand{\epr}{\hspace{\fill}$\Box$}

\newcommand{\bg}{\bigskip}
\newcommand{\med}{\medskip}
\newcommand{\sm}{\smallskip}

\font\caps=cmcsc10 scaled \magstep1   

\begin{document}
\noindent {\em Semigroup Forum} (2011) {\bf 83}, 250 -- 280\\
DOI 10.1007/s00233-011-9304-z
\vspace{0.02in}\\
{\em arXiv version}: layout, fonts, pagination and numbering of sections, lemmas, theorems and formulas may vary from the SF published version
\vspace{0.04in}\\
\title[Lattice isomorphisms of bisimple monogenic orthodox semigroups]{Lattice isomorphisms of bisimple\\ monogenic orthodox semigroups}
\author{Simon M. Goberstein}
\address{\hspace*{-0.205in}Department of Mathematics \!\&\! Statistics,
California State University, Chico, CA 95929, USA, \!\!e-mail: sgoberstein@csuchico.edu}
\begin{abstract}
Using the classification and description of the structure of bisimple monogenic orthodox semigroups obtained in \cite{key10}, we prove that every bisimple orthodox semigroup generated by a pair of mutually inverse elements of infinite order is strongly determined by the lattice of its subsemigroups in the class of all semigroups. This theorem substantially extends an earlier result of \cite{key25} stating that the bicyclic semigroup is strongly lattice determined. 
\vspace{0.07in}
\end{abstract}
\maketitle
\font\caps=cmcsc10 scaled \magstep1   
\def\bfseries{\normalsize\caps}
\vspace{-0.28in}
\section*{Introduction}
\sm
Let $S$ be a semigroup. The set of all subsemigroups of $S$ (including, by convention, the empty one) is a lattice under set-theoretic inclusion, and the relationship between the properties of this lattice and the properties of $S$ has been studied in numerous publications for over half a century. One of the central problems in this area of research is deciding whether $S$ is isomorphic or antiisomorphic to any semigroup $T$ whose lattice of subsemigroups is isomorphic to that of $S$, and if this is the case, $S$ is said to be lattice determined in the class of all semigroups. If for each isomorphism $\Phi$ of the subsemigroup lattice of $S$ onto that of a semigroup $T$, there exists an isomorphism or an antiisomorphism $\varphi$ of $S$ onto $T$ such that $U\Phi=U\varphi$ for every subsemigroup $U$ of $S$, then $S$ is called strongly lattice determined. For instance, it has long been known that the infinite cyclic group is strongly lattice determined \cite[Lemma 34.8]{key30}. We refer to \cite{key30} for a comprehensive treatment of results concerning lattice determinability of semigroups of various types obtained prior to 1996. As the content of \cite{key30} shows, the problem of lattice determinability of a {\em regular} semigroup $S$ in the class of {\em all} semigroups had been considered {\em systematically} only in the cases when $S$ is a group, a semilattice, a rectangular band or, more generally, a completely simple semigroup \cite[Sections 34, 36, 37, and 38, respectively]{key30}; more recently, definitive results about lattice determinability of completely $0$-simple semigroups were obtained in \cite{key20}. In addition, it was proved in \cite{key25} (see also \cite[Theorem 41.8]{key30}) that the bicyclic semigroup is strongly lattice determined. The main goal of this paper is to extend the latter result to the class of all bisimple orthodox semigroups generated by a pair of mutually inverse nongroup elements.

By analogy with the inverse semigroup case, we call an orthodox semigroup monogenic if it is generated by two mutually inverse elements. It is well known that cyclic groups and the bicyclic semigroup are the only bisimple monogenic inverse semigroups. However, the class of bisimple monogenic orthodox semigroups is substantially more diverse. A complete classification and description of the structure of semigroups of that class was obtained only recently by the author \cite{key10}. In particular, we constructed in \cite{key10} a family of pairwise nonisomorphic bisimple orthodox semigroups $\cO_{(\nu,\,\mu)}(a,b)$ indexed by ordered pairs $(\nu,\mu)\in\N^{\ast}\times\N^{\ast}$ (here and elsewhere in this paper $\N$ denotes the set of all positive integers and $\N^{\ast}=\N\cup\{\infty\}$), each being generated by a pair of mutually inverse elements $a$ and $b$ satisfying $ab=a^2b^2$ and $ba\neq b^2a^2$, and proved that if $S$ is an arbitrary bisimple monogenic orthodox semigroup with nongroup generators, then $S$ or its dual is isomorphic to one of the semigroups of that two-parameter family. Needless to say that the results of \cite{key10} play a crucial role in this paper. 

The paper is organized as follows. Section 1 contains basic information about orthodox semigroups and lattice isomorphisms of semigroups plus a few new auxiliary results some of which (for instance, Proposition \ref{108}) might be of independent interest. In Section 2 we review the results of \cite{key10} about bisimple monogenic orthodox semigroups with nongroup generators and establish a number of additional useful properties of these semigroups. In Section 3 we prove our principal new theorem: every bisimple orthodox semigroup generated by a pair of mutually inverse elements of infinite order is strongly lattice determined (Theorem \ref{301}). 

The main results of the paper were reported at the International Conference on Geometric and Combinatorial Methods in Group Theory and Semigroup Theory held at the University of Nebraska, Lincoln, on May 17-21, 2009.    
\sm
\section{Preliminaries}
\sm
We use the term ``order'' instead of ``partial order'' and refer to a linearly ordered set as a chain. If $(A,\leq)$ is an ordered set and $a\in A$, the principal order ideal of $A$ generated by $a$ will be denoted by $(a\;\!]$, so $(a\;\!]=\{b\in A:b\leq a\}$.  Let $S$ be a semigroup. We say that $x\in S$ is a {\em group element} of $S$ if it belongs to some subgroup of $S$; otherwise $x$ is a {\em nongroup element} of $S$. The set of nongroup elements of $S$ will be denoted by $N_S$, and the set of idempotents of $S$ by $E_S$. To indicate that $U$ is a subsemigroup of $S$, we will write $U\leq S$. Under the convention that $\emptyset\leq S$, the set of all subsemigroups of $S$ ordered by inclusion is a lattice which we denote by $\cSb(S)$. Clearly, $\cSb(S)=\cSb(S^{\rm opp})$ where $S^{\rm opp}$ is the dual (or ``opposite'') semigroup of $S$. As usual, $\langle X\rangle$ stands for the subsemigroup of $S$ generated by $X\subseteq S$, and $\langle x\rangle$ for the cyclic subsemigroup of $S$ generated by $x\in S$.  The order of an element $x$ of $S$ is denoted by $o(x)$; if $x$ has infinite order, we write $o(x)=\infty$. If $w=w(x_1,\ldots,x_n)$ is a word in the alphabet $\{x_1,\ldots,x_n\}\subseteq S$, we will say that $w$ is a word in $x_1,\ldots, x_n$ and will identify $w$ with its value in $S$ if no confusion is likely to occur. For any $x\in S$, we denote by $x^0$ the identity of $S^1$, so $x^0y=yx^0=y$ for all $y\in S$ (exceptions of this agreement might happen when $S$ contains a subsemigroup $U$ with an identity $e$ and we put $x^0=e$ for each $x\in U$; all such situations will be clear from the context). As in \cite{key10}, to indicate that $x,y\in S$ satisfy $xyx=x$ and $yxy=y$, we will write $x\perp y$, and the phrase ``$x\perp y$ in $S$'' will mean that $x,y\in S$ and $x\perp y$. We also denote by $V_S(x)$ the set of all inverses of $x\in S$, so $y\in V_S(x)$ if and only if $x\perp y$ in $S$. Standard facts about Green's relations on semigroups will be used without reference. If $U\leq S$, we distinguish Green's relations on $U$ from those on $S$ by using superscripts. According to \cite[Result 9]{key14}, if $U$ is a regular subsemigroup of $S$, then $\cK^U=\cK^S\cap(U\times U)$ for $\cK\in\{\cL, \cR, \cH\}$. This result will also be applied without mention. It is common to say that $S$ is {\em combinatorial} if $\cH$ is the identity relation on $S$, so a regular semigroup is combinatorial if and only if it has no nontrivial subgroups.

Recall that an {\em orthodox semigroup} is a regular semigroup in which the idempotents form a subsemigroup. By \cite[Theorem VI.1.1]{key15}, if $S$ is a regular semigroup, the following conditions are equivalent: (a) $S$ is orthodox; (b) $V_S(e)\subseteq E_S$ for all $e\in E_S$; (c) $V_S(b)V_S(a)\subseteq V_S(ab)$ for all $a,b\in S$. Thus if $x\perp y$ in an orthodox semigroup $S$, then $x\in N_S$ if and only if $y\in N_S$, and $x^n\perp y^n$ for all $n\in\N$, which implies that $o(x)=o(y)$. These simple facts will be used without comment. According to the terminology introduced in \cite{key10}, $S$ is a {\em monogenic orthodox semigroup} if it is an orthodox semigroup generated by a pair of mutually inverse elements. In what follows, the phrase ``let $S=\langle a,b\rangle$ be a monogenic orthodox semigroup'' will always mean that $S$ is an orthodox semigroup with $a\perp b$ in $S$.

As in \cite{key30}, left [right] zero semigroups will be called {\em left {\rm[}right\,{\rm]} singular}, and a semigroup is {\em singular} if it is either left or right singular. By definition \cite[\S 1.8]{key5}, a {\em rectangular band} is the Cartesian product $I\times\Lambda$ where $I$ is a left and $\Lambda$ a right singular semigroup. If $S$ is a semigroup such that $S=\bigcup_{(i,\lambda)\in I\times\Lambda} S_{i\lambda}$ where $S_{i\lambda}\leq S$ for all $(i,\lambda)\in I\times\Lambda$, and $S_{i\lambda}\cap S_{j\mu}=\emptyset$ whenever $(i,\lambda)\not=(j,\mu)$, then $S$ is a {\em rectangular band of semigroups $S_{i\lambda}$} provided that $S_{i\lambda}S_{j\mu}\subseteq S_{i\mu}$ for all $(i,\lambda),(j,\mu)\in I\times\Lambda$, and if the set $I$ [$\Lambda$] is a singleton, $S$ is also called a {\em right} [{\em left}\,] {\em singular  band of $S_{i\lambda}$} (see \cite[\S 1]{key30}).

\indent The following technical fact is probably well-known. We record it for convenience of reference and include its proof for completeness.
 
\bl\label{101} Let $S$ be a semigroup, and let $x,y\in S$. If $xy=y^h$ {\rm[}$xy=x^h${\rm]} for some $h\in\N$, then $x^my^n=y^{m(h-1)+n}$ {\rm[}$x^my^n=x^{n(h-1)+m}${\rm]} for all $m,n\in\N$.
\el
{\bf Proof.} Suppose $xy=y^h$ for some $h\in\N$. Take any $n\in\N$. Then $xy^n=(xy)y^{n-1}=y^{h-1+n}$, and for any $k>1$ such that $x^{k-1}y^n=y^{(k-1)(h-1)+n}$, we have 
\[x^ky^n=x(x^{k-1}y^n)=xy^{(k-1)(h-1)+n}=(xy)y^{k(h-1)-h+n}=y^hy^{k(h-1)-h+n}=y^{k(h-1)+n},\]
so $x^my^n=y^{m(h-1)+n}$ for all $m\in\N$. We have shown that $x^my^n=y^{m(h-1)+n}$ for all $m,n\in\N$, and the alternative statement holds by duality. \epr

\bl\label{102} Let $S$ be a semigroup, and let $x$ and $y$ be two distinct elements of $S$ such that $o(x)=o(y)=\infty$, $xy=y^2$, and $yx=x^2$. Then
\vspace{0.03in}\\
\indent{\rm(i)} $x^my^n=y^{m+n}$ and $y^nx^m=x^{n+m}$ for all $m,n\in\N$;
\vspace{0.02in}\\
\indent{\rm(ii)} $\langle x,y\rangle=\langle x\rangle\cup\langle y\rangle$;
\vspace{0.02in}\\
\indent{\rm(iii)} either $\langle x\rangle\cap\langle y\rangle=\emptyset$, in which case $\langle x,y\rangle$ is a right singular band of $\langle x\rangle$ and $\langle y\rangle$, or there is $m\geq 2$ such that $x^m=y^m$ but $x^i\neq y^i$ if $1\leq i\leq m-1$, in which case $x^n=y^n$ for all $n\geq m$ and $x^k\neq y^{l}$ for all $k,l\in\N$ such that $k\neq l$.
\el
{\bf Proof.} (i) This is an immediate corollary of Lemma \ref{101}. 
\vspace{0.02in}\\
\indent(ii) Let $z\in\langle x,y\rangle$. Then $z=x^{k_1}_1\cdots x^{k_r}_r$ for some $r\geq 1$ and $k_1,\ldots,k_r\in\N$ where $x_i\in\{x,y\}$ for $1\leq i\leq r$ and $x_j\neq x_{j+1}$ for $1\leq j\leq r-1$. By induction on $r$, it is easily seen that $z=x^{k_1+\cdots+k_r}_r\in\langle x\rangle\cup\langle y\rangle$. Therefore $\langle x,y\rangle=\langle x\rangle\cup\langle y\rangle$.
\vspace{0.02in}\\
\indent(iii) If $\langle x\rangle\cap\langle y\rangle=\emptyset$, it follows from (i) and (ii) that $\langle x,y\rangle$ is a right singular band of $\langle x\rangle$ and $\langle y\rangle$. Suppose that $\langle x\rangle\cap\langle y\rangle\neq\emptyset$. Then $x^k=y^l$ for some $k,{l}\in\N$ whence, using (i), we obtain $y^{k+l}=x^ky^l=y^{2l}$, so that $k=l$ since $o(y)=\infty$. Thus $x^k\neq y^{l}$ if $k\neq {l}$ in $\N$. Let $m=\min\{k\in\N:x^k=y^k\}$. Since $x\neq y$, it is clear that $m\geq 2$ and $x^i\neq y^i$ if $1\leq i\leq m-1$. Using (i) again, for each $n\geq m$, we have $x^n=x^{n-m}x^m=x^{n-m}y^m=y^n$.\epr 
\vspace{0.05in}\\
\indent If $S$ and $T$ are semigroups such that $\cSb(S)\cong\cSb(T)$, then $S$ and $T$ are called {\em lattice isomorphic}, and any isomorphism of $\cSb(S)$ onto $\cSb(T)$ is referred to as a {\em lattice isomorphism} of $S$ onto $T$. If $\Phi$ is a lattice isomorphism of $S$ onto $T$, then $\Phi$ is said to be {\em induced} by a mapping $\varphi : S\rightarrow T$ if $U\Phi=U\varphi$ for all $U\leq S$. If $S$ is isomorphic or antiisomorphic to any semigroup that is lattice isomorphic to it, then $S$ is called {\em lattice determined}, and if each lattice isomorphism of $S$ onto a semigroup $T$ is induced by an isomorphism or an antiisomorphism of $S$ onto $T$, we say that $S$ is {\em strongly lattice determined}.
 
\br\label{103}{\rm (From \cite[\S 31]{key30})} Let $S$ be a semigroup in which every nonidempotent element has infinite order, and let $\Phi$ be a lattice isomorphism of $S$ onto a semigroup $T$. Then $\Phi$ is induced by a unique bijection $\varphi : S\rightarrow T$ which is defined by the formula $\langle x\rangle\Phi=\langle x\varphi\rangle$ for all $x\in S$ and has the property that $(x^n)\varphi=(x\varphi)^n$ for all $x\in S$ and $n\in\N$.
\er 
If $S$ is a semigroup in which every nonidempotent element has infinite order and if $\Phi$ is a lattice isomorphism of $S$ onto a semigroup $T$, then the bijection $\varphi : S\rightarrow T$ described in Result \ref{103} will be called the {\em $\Phi$-associated bijection} of $S$ onto $T$.

The {\em bicyclic semigroup} $\cB(a,b)$ is usually defined as a semigroup with identity $1$ generated by the two-element set $\{a,b\}$ and given by one defining relation $ab=1$ (see \cite[\S 1.12]{key5}). Clearly, $\cB(a,b)$ can also be defined without mentioning the identity as a semigroup given by the following presentation: $\cB(a,b)=\langle a, b\;|\;aba=a, bab=b, a^2b=a, ab^2=b\rangle$. The structure of $\cB(a,b)$ is well known -- it is a combinatorial bisimple inverse semigroup, each of its elements has a unique representation in the form $b^ma^n$ where $m$ and $n$ are nonnegative integers (and $a^0=b^0=ab$), the semilattice of idempotents of $\cB(a,b)$ is a chain: $ab>ba>b^2a^2>\cdots$, and for a given $b^ma^n\in\cB(a,b)$, we have $R_{b^ma^n}=\{b^ma^l\,:\,l\geq 0\}$ and $L_{b^ma^n}=\{b^ka^n\,:\,k\geq 0\}$ (see \cite[Lemma 1.31 and Theorem 2.53]{key5}). Since each nonidempotent element of $\cB(a,b)$ has infinite order, by Result \ref{103} any lattice isomorphism $\Phi$ of $\cB(a,b)$ onto a semigroup $T$ is induced by the $\Phi$-associated bijection $\varphi$, and by \cite[Main Theorem]{key25} (or \cite[Theorem 41.8]{key30}), $\varphi$ is an isomorphism or an antiisomorphism of $\cB(a,b)$ onto $T$. In other words, we have 
 
\br\label{104}{\rm\cite[Main Theorem]{key25}} The bicyclic semigroup is strongly lattice determined.
\er

It is immediate from \cite[Proposition 1]{key23} (reproduced as \cite[Proposition 3.2(a)]{key30}) that a semigroup which is lattice isomorphic to a band is itself a band. On the other hand, it is also easily seen (and well known) that bands, in general, are not lattice determined -- for instance, it is clear that any singular semigroup $S$ is lattice isomorphic to a chain having the same cardinality as $S$. However, according to \cite[Theorem 7]{key23} (see also \cite[Theorem 37.8]{key30}), the situation is quite different for nonsingular rectangular bands.

\br\label{105}{\rm\cite[Theorem 7]{key23}} Any nonsingular rectangular band is strongly lattice determined.
\er

Note that {\em any} rectangular band of nonperiodic abelian groups is strongly lattice determined (even if that band is singular). In fact, a more general theorem holds.

\br\label{106}{\rm \cite[Main Theorem]{key1}} A rectangular band of nonperiodic commutative cancellative semigroups without adjoined identities is strongly lattice determined.
\er

In Lemma \ref{107} below, we consider lattice isomorphisms of a semigroup $T$ generated by two distinct elements $s$ and $t$ of infinite order such that $st=t^2$ and $ts=s^2$ (clearly, $T$ is commutative precisely when $s^2=t^2$). Part of the proof of this lemma draws on techniques used in the proof of \cite[Lemma 35.4]{key30}, which describes properties of the $\Phi$-associated bijection where $\Phi$ is a lattice isomorphism of a {\em commutative} idempotent-free semigroup, in which no element is an identity for another element and which contains elements $x$ and $y$ satisfying $xy=x^m$ for some $m\geq 2$, onto some {\em commutative} semigroup. Although commutativity of semigroups is essentially used in the proof of \cite[Lemma 35.4]{key30}, several of its technical arguments work in a noncommutative situation, and we use them in the proof of our lemma (see the proof of \cite[Lemma 35.4]{key30} for comparison).
 
\bl\label{107} Let $T=\langle s,t\rangle$ be a semigroup such that $s\neq t$, $o(s)=o(t)=\infty$, $st=t^2$, and $ts=s^2$. Let $\Phi$ be a lattice isomorphism of $T$ onto a semigroup $Q$, and let $\varphi$ be the $\Phi$-associated bijection of $T$ onto $Q$. Then one of the following holds:
\vspace{0.02in}\\
\indent{\rm(i)} $\varphi$ is an isomorphism or an antiisomorphism of $T$ onto $Q$, or
\vspace{0.02in}\\ 
\indent{\rm(ii)} $s^2\neq t^2$ but $s^n=t^n$ for all integers $n\geq 3$, and either $(s\varphi)(t\varphi)=(t\varphi)(s\varphi)=(t\varphi)^2$ or $(s\varphi)(t\varphi)=(t\varphi)(s\varphi)=(s\varphi)^2$.
\el
{\bf Proof.} Denote $p=s\varphi$ and $q=t\varphi$. By Result \ref{103}, $o(p)=o(q)=\infty$. By Lemma \ref{102}(ii), $T=\langle s\rangle\cup\langle t\rangle$ and therefore, according to \cite[Lemma 31.7]{key30}, $Q=\langle p\rangle\cup\langle q\rangle$. By Lemma \ref{102}(iii), if $\langle s\rangle\cap\langle t\rangle=\emptyset$, then $T$ is a right singular band of $\langle s\rangle$ and $\langle t\rangle$, so that, by Result \ref{106}, $\varphi$ is an isomorphism or an antiisomorphism of $T$ onto $Q$. Now suppose that $\langle s\rangle\cap\langle t\rangle\neq\emptyset$. Then, by Lemma \ref{102}(iii), there is $m\geq 2$ such that $s^n=t^n$ for all $n\geq m$ but $s^k\neq t^{l}$ for all $k,l\in\N$ such that either $k\neq l$ or $k=l\in\{1,\ldots, m-1\}$. By Result \ref{103}, it follows that $p^n=q^n$ for all $n\geq m$ and $p^k\neq q^{l}$ if $k\neq l$ or if $k=l\in\{1,\ldots, m-1\}$. If $pq=p$, then using Lemma \ref{101}, we obtain $p^{m+1}=p\cdot p^m=pq^m=p^{m(1-1)+1}=p$, which contradicts the fact that $o(p)=\infty$. Therefore $pq\neq p$ and, dually, $qp\neq p$. By symmetry, we also have $qp\neq q$ and $pq\neq q$. Thus $pq, qp\in\{p^2,p^3,\ldots\}\cup\{q^2,q^3,\ldots\}$.

Suppose $pq\in\{q^2,q^3,\ldots\}$. Then $pq=q^k$ for some $k\geq 2$, and hence, by Lemma \ref{101}, $q^{m+1}=p^{m-1}(pq)=p^{m-1}q^k=q^{(m-1)(k-1)+k}=q^{mk-m+1}$. Since $o(q)=\infty$, it follows that $k=2$, that is, $pq=q^2$. Thus $pq\in\{q^2,q^3,\ldots\}$ if and only if $pq=q^2$. Dually, $qp\in\{q^2,q^3,\ldots\}$ if and only if $qp=q^2$. Therefore, by symmetry, $qp\in\{p^2,p^3,\ldots\}$ if and only if $qp=p^2$, and $pq\in\{p^2,p^3,\ldots\}$ if and only if $pq=p^2$. 
\vspace{0.02in}\\
\indent {\bf Case 1:} $pq=q^2$.
\vspace{0.01in}\\
\indent We consider separately the following situations: (a) $m=2$, (b) $m\geq 4$, and (c) $m=3$.
\vspace{0.02in}\\
\indent 1(a) $m=2$.
\vspace{0.01in}\\
\indent In this case, $p^n=q^n$ for all $n\geq 2$. Hence $pq=q^2=p^2$ and, dually, $qp=p^2=q^2$, from which it is immediate that $\varphi$ is an isomorphism of $T$ onto $Q$.
\vspace{0.02in}\\
\indent 1(b) $m\geq 4$.
\vspace{0.01in}\\
\indent As noted above, either $qp=p^2$ or $qp=q^2$. Let us show that the latter cannot happen here. Suppose $qp=q^2$. Then $pq=q^2=qp$ and so $Q$ is commutative. In particular, $p^2q=qp^2$. Let $r\in\langle p^2,q\rangle$. Then $r=p^{2i}q^j$ for some integers $i,j\geq 0$ such that $i+j\geq 1$. If $j=0$, then $r=p^{2i}\in\langle p^2\rangle$. Now suppose that $j\geq 1$. If $i=0$, then $r=q^j\in\langle q\rangle$, and if $i\geq 1$, applying Lemma \ref{102}(i), we obtain $r=p^{2i}q^j=q^{2i+j}\in\langle q\rangle$. Therefore $\langle p^2,q\rangle=\langle p^2\rangle\cup\langle q\rangle$. Since $\langle p^2,q\rangle\Phi^{-1}=\langle s^2,t\rangle$, according to \cite[Lemma 31.7]{key30}, $\langle s^2,t\rangle=\langle p^2\rangle\Phi^{-1}\cup\langle q\rangle\Phi^{-1}=\langle s^2\rangle\cup\langle t\rangle$. By Lemma \ref{102}(i), $s^3=ts^2\in\langle s^2,t\rangle=\langle s^2\rangle\cup\langle t\rangle$. Since $s^3\not\in\langle s^2\rangle$, we have $s^3\in\langle t\rangle$, contradicting the assumption that $m\geq 4$ and the fact that $s^k\neq t^{l}$ if $k\neq l$ or $k=l\in\{1,\ldots,m-1\}$.

Since $qp\neq q^2$, we must have $qp=p^2$. By Lemma \ref{102}(i), $p^kq^{l}=q^{k+l}$ and $q^{l}p^k=p^{l+k}$ for all $k,l\in\N$. Take any $x,y\in T=\langle s\rangle\cup\langle t\rangle$. If $x,y\in\langle s\rangle$ or $x,y\in\langle t\rangle$, it is clear that $(xy)\varphi=(x\varphi)(y\varphi)$. Suppose that $x\in\langle s\rangle$ and $y\in\langle t\rangle$. Then $x=s^k$ and $y=t^{l}$ for some $k,l\in\N$. Hence $(xy)\varphi=(s^kt^{l})\varphi=(t^{k+l})\varphi=q^{k+l}=p^kq^{l}=(x\varphi)(y\varphi)$ and, similarly, $(yx)\varphi\!=\!(t^{l}s^k)\varphi\!=\!(s^{l+k})\varphi\!=\!p^{l+k}\!=\!q^{l}p^k\!=\!(y\varphi)(x\varphi)$. Thus $\varphi$ is an isomorphism of $T$ onto $Q$.
\vspace{0.02in}\\
\indent 1(c) $m=3$.
\vspace{0.01in}\\
\indent Here $p^n=q^n$ for all $n\geq 3$, $\{p,p^2\}\cap\{q,q^2\}\!=\!\emptyset$, and once again either $qp\!=\!p^2$ or $qp\!=\!q^2$. If $qp=p^2$, then exactly as in 1(b) we can show that $\varphi$ is an isomorphism of $T$ onto $Q$. Now let $qp=q^2$. Then $pq=q^2=qp$, that is, $(s\varphi)(t\varphi)=(t\varphi)(s\varphi)=(t\varphi)^2$. Unlike a similar situation considered in 1(b), here the assumption that $qp=q^2$ does not lead to a contradiction. Note that in this case lattice isomorphic semigroups $T$ and $Q$ are neither isomorphic nor antiisomorphic since $Q$ is commutative but $T$ is not.
\vspace{0.02in}\\
\indent {\bf Case 2:} $pq=p^2$.
\vspace{0.02in}\\
\indent If $m=2$, then exactly as in 1(a), we observe that $\varphi$ is an isomorphism of $T$ onto $Q$. If $m\geq 4$, by the dual of the argument used in 1(b), $\varphi$ is an antiisomorphism of $T$ onto $Q$. Finally, if $m=3$, then dually to 1(c), either $\varphi$ is an antiisomorphism of $T$ onto $Q$, or $Q$ is commutative and $(s\varphi)(t\varphi)=(t\varphi)(s\varphi)=(s\varphi)^2$. \epr
\bp\label{108} Let $S=\langle s,f\rangle$ be a semigroup such that $f^2=f$, $o(s)=\infty$, $fs=s$, and $sf\neq s$. Let $\Phi$ be a lattice isomorphism of $S$ onto a semigroup $P$, and let $\varphi$ be the $\Phi$-associated bijection of $S$ onto $P$. Then $f\varphi$ is an adjoined identity in the subsemigroup $\langle f\varphi, (sf)\varphi\rangle$ of $P$ and one of the following two statements is true:
\vspace{0.03in}\\
\indent{\rm(I)} $\varphi$ is an isomorphism or an antiisomorphism of $S$ onto $P$, or
\vspace{0.03in}\\
\indent{\rm(II)} $(s\varphi)(f\varphi)=(f\varphi)(s\varphi)=(sf)\varphi$,
\vspace{0.03in}\\
with {\rm(II)} being possible only if one of the following holds: either $s^nf=s^n$ for all $n\geq 2$, which implies $(s\varphi)\cdot(sf)\varphi=(sf)\varphi\cdot(s\varphi)=[(sf)\varphi]^2=(s\varphi)^2$ and so $[(sf)\varphi]^n=(s\varphi)^n$ for all $n\geq 2$, or $s^2f\neq s^2$ but $s^nf=s^n$ for all $n\geq 3$, in which case $(s\varphi)\cdot(sf)\varphi=(sf)\varphi\cdot(s\varphi)=[(sf)\varphi]^2\neq(s\varphi)^2$ and $[(sf)\varphi]^n=(s\varphi)^n$ for all $n\geq 3$.
\ep
{\bf Proof.} Let $t=sf$ and $T=\langle s,t\rangle$. Then $tf=t$ and, by assumption, $t\neq s$ and $fs=s$. It follows that $ft=t$, $t^n=(sf)^n=s^nf$ for all $n\in\N$, and $f\notin\langle s\rangle$ since $o(s)=\infty$. In particular, we have $st=s^2f=t^2$ and $ts=sfs=s^2$. If $t^k=t^{l}$ for some $k,l\in\N$, then $s^{k+1}=t^ks=t^{l}s=s^{l+1}$ and so $k=l$ since $o(s)=\infty$. Thus $o(t)=\infty$ and $f$ is an adjoined identity in the subsemigroup $\langle f,t\rangle$ of $S$; in particular, $\langle f,t\rangle=\{f\}\cup\langle t\rangle$. By Lemma \ref{102}(ii), $T=\langle s\rangle\cup\langle t\rangle$, so that $S=\{f\}\cup\langle s\rangle\cup\langle t\rangle$ and $f\notin\langle s\rangle\cup\langle t\rangle$.

Let $e=f\varphi$, $p=s\varphi$, $q=t\varphi$, and $Q=T\Phi$. Then $e^2=e$, $o(p)=o(q)=\infty$, $Q=\langle p,q\rangle$, and $P=\langle p,e\rangle$. By \cite[Lemma 31.7]{key30}, $Q=\langle p\rangle\cup\langle q\rangle$. Therefore $P=\{e\}\cup\langle p\rangle\cup\langle q\rangle$, $e\notin\langle p\rangle\cup\langle q\rangle$, and either $pe\in\langle q\rangle$ or $ep\in\langle q\rangle$. Since $\Phi|_{\cSb(\langle f,t\rangle)}$ is a lattice isomorphism of $\langle f, t\rangle$ onto $\langle e, q\rangle$, it follows from \cite[Lemma 33.4]{key30} that $e$ is either an adjoined identity or an adjoined zero in $\langle e, q\rangle$. The latter leads to a contradiction since $pe\in\langle q\rangle$ implies $pe=(pe)e=e\not\in\langle q\rangle$ and, dually, $ep\in\langle q\rangle$ implies $ep\not\in\langle q\rangle$. Hence $e$ is an adjoined identity in $\langle e, q\rangle$.

Suppose that $pq=q^2$. If $ep=e$, then $q^2=eq^2=(ep)q=eq=q$, and if $pe=e$, then $q^2=pq=(pe)q=eq=q$, contrary to the fact that $o(q)=\infty$. Thus
\setcounter{equation}{8} 
\be
pq=q^2\Longrightarrow ep\neq e\text{ and }pe\neq e.
\label{eq:1st}
\ee
If $pe=q^k$ for some $k\geq 1$, then $q^{k+1}=peq=pq=q^2$ whence $k=1$ and therefore $pe=q$. Similarly, if $ep=q^k$ for some $k\geq 1$, then $q^{k+1}=epq=eq^2=q^2$, so that $k=1$ and $ep=q$. Since either $pe\in\langle q\rangle$ or $ep\in\langle q\rangle$, we have shown that 
\be 
pq=q^2\Longrightarrow pe=q\text{ or }ep=q.
\label{eq:2nd}
\ee
 
Note that if $qp=p^2$ and $ep=p^k$ for some $k\!\in\!\N$, then $p^2=qep=qp^k=(qp)p^{k-1}=p^{k+1}$ whence $k=1$ and so $ep=p$. Therefore 
\be
qp=p^2\text{ and }ep\in\langle p\rangle\Longrightarrow ep=p.
\label{eq:3rd}
\ee

By Lemma \ref{107}, since $\Phi|_{\cSb(T)}$ is a lattice isomorphism of $T$ onto $Q$, one of the following is true: (i) $\varphi|_T$ is an isomorphism or an antiisomorphism of $T$ onto $Q$, or (ii) $s^2\neq t^2$ but $s^n=t^n$ for all $n\geq 3$, and either $(s\varphi)(t\varphi)=(t\varphi)(s\varphi)=(t\varphi)^2$ or $(s\varphi)(t\varphi)=(t\varphi)(s\varphi)=(s\varphi)^2$. 

Suppose that (ii) holds, that is, $s^2\neq s^2f$ whereas $s^n=s^nf$ for all $n\geq 3$, and either $pq=qp=q^2$ or $pq=qp=p^2$. Assume that $pq=qp=p^2$. Then for all $m\geq 2$, we have $p^me=p^{m-2}(pqe)=p^{m-2}(pq)=p^m$ and, dually, $ep^m=p^m$, so $e$ is an adjoined identity in the semigroup $\{e\}\cup\{p^2,p^3,\ldots\}$. Since $fs^2=s^2\neq f$, by \cite[Lemma 33.4]{key30} $f$ is an adjoined identity in the semigroup $\{f\}\cup\{s^2,s^3,\ldots\}$ whence $s^2=s^2f$; a contradiction. Thus $pq=qp=q^2$, that is, $(s\varphi)\cdot(sf)\varphi=(sf)\varphi\cdot(s\varphi)=[(sf)\varphi]^2\neq(s\varphi)^2$. By (\ref{eq:2nd}), $pe=q$ or $ep=q$. Suppose $pe=q$ but $ep\neq q$. As noted prior to (\ref{eq:2nd}), the latter implies $ep\notin\langle q\rangle$ and so, in view of (\ref{eq:1st}), $ep=p^k$ for some $k\in\N$. Hence $p^{k+1}=pep=qp=q^2$. If $k\geq 2$, then $q^{k+1}=p^{k+1}=q^2$; a contradiction. Thus $k=1$ and $ep=p$. Then $p^2=p(ep)=(pe)p=qp=q^2$, again a contradiction. Therefore $pe=q$ and $ep\neq q$ cannot hold simultaneously. By duality, $ep=q$ and $pe\neq q$ cannot hold simultaneously as well. It follows that $pe=ep=q$, that is, $(s\varphi)(f\varphi)=(f\varphi)(s\varphi)=(sf)\varphi$.

For the remainder of the proof, we will assume that (i) holds. 

Suppose that $s^2=s^2f$. Then $pq=qp=p^2=q^2$, so $T$ and $Q$ are commutative semigroups and $\varphi|_{T}$ is an isomorphism of $T$ onto $Q$. By (\ref{eq:2nd}), $pe=q$ or $ep=q$. If $pe=q$ but $ep\neq q$, by (\ref{eq:1st}) and (\ref{eq:3rd}) we have $ep=p$ whence $(s\varphi)(f\varphi)=(sf)\varphi$ and $(f\varphi)(s\varphi)=s\varphi=(fs)\varphi$, and since $\varphi|_{T}$ is an isomorphism of $T$ onto $Q$, it follows that $\varphi$ is an isomorphism of $S$ onto $P$. If $ep=q$ but $pe\neq q$, then $pe=p$ by (\ref{eq:1st}) and the dual of (\ref{eq:3rd}), so that $(f\varphi)(s\varphi)=(sf)\varphi$ and $(s\varphi)(f\varphi)=s\varphi=(fs)\varphi$, and therefore $\varphi$ is an antiisomorphism of $S$ onto $P$. On the other hand, if $pe=ep=q$, then $(s\varphi)(f\varphi)=(f\varphi)(s\varphi)=(sf)\varphi$ and we conclude that $\varphi$ is neither an isomorphism nor an antiisomorphism of $S$ onto $P$. 

Now assume that $s^2\neq s^2f$ and hence $p^2\neq q^2$.
\vspace{0.02in}\\
\indent {\bf Case 1:} $\varphi|_T$ is an isomorphism of $T$ onto $Q$.
\vspace{0.02in}\\
\indent Here $pq=(s\varphi)(t\varphi)=(st)\varphi=(t^2)\varphi=q^2$ and $qp=(t\varphi)(s\varphi)=(ts)\varphi=(s^2)\varphi=p^2$. As observed prior to (\ref{eq:2nd}), if $ep\in\langle q\rangle$, then $ep=q$, which implies $p^2=qp=(qe)p=q(ep)=q^2$; a contradiction. Thus $ep\not\in\langle q\rangle$. By (\ref{eq:1st}), $ep\neq e$ and so $ep\in\langle p\rangle$. In view of (\ref{eq:3rd}), $ep=p$ and hence, by (\ref{eq:2nd}), $pe=q$. Therefore $(f\varphi)(s\varphi)=s\varphi=(fs)\varphi$ and $(s\varphi)(f\varphi)=(sf)\varphi$. Since $\varphi|_T$ is an isomorphism of $T$ onto $Q$, it follows that $\varphi$ is an isomorphism of $S$ onto $P$. 
\vspace{0.02in}\\
\indent {\bf Case 2:} $\varphi|_T$ is an antiisomorphism of $T$ onto $Q$. 
\vspace{0.02in}\\
\indent By the argument dual to that of Case 1, $\varphi$ is an antiisomorphism of $S$ onto $P$.  \epr
\sm
\section{Bisimple monogenic orthodox semigroups with nongroup generators}
\sm
Let $S$ be an orthodox semigroup, $a\in N_S$, $b\in V(a)$, and $ab=a^2b^2$. If $m,n\in\Z$ and $i,j\in\{0,1\}$ are such that either (I) $m>i$ and $n=j=0$, or (II) $m=i=0$ and $n\geq 1$, or (III) $m>i$ and $n>j$, then in the terminology of \cite{key10} $a^ib^ma^nb^j$ is an {\em abridged word in $a, b$} (in this order!) of {\em type I, II, or III}, respectively (or simply an {\em abridged word} in $a, b$ when there is no need to indicate its type); in order to avoid repetition, we sometimes omit the phrase ``in $a,b$'' if no confusion is likely. In this paper, instead of calling $a^ib^ma^nb^j$ an abridged word in $a,b$ (in this order), we will often refer to it as an {\em abridged (a,b)-word} (of a particular type). We will say that $x\in S$ is {\em represented} by an abridged $(a,b)$-word $a^ib^ma^nb^j$ if $x=a^ib^ma^nb^j$, in which case we will also refer to $a^ib^ma^nb^j$ as an {\em abridged form} of $x$, specifying its type if necessary. As noted in \cite[Remark 1]{key10}, the definitions of type I and type II abridged words in $a,b$ are not entirely symmetric -- in a type I word $a^ib^m$ we must have $m>i$ but in a type II word $a^nb^j$ the equality $n=j=1$ is allowed. This ensures that if $x=ab$, there is a unique abridged $(a,b)$-word representing $x$, namely, a type II word $a^0b^0a^1b^1$. Thus the dual of an abridged $(b,a)$-word $b^1a^1$ of type II is an abridged $(a,b)$-word $a^1b^1$ of type II; this is the only exception of the following rule: the dual of an abridged $(b,a)$-word of type I [II, III] is an abridged $(a,b)$-word of type II [I, III]. We will call $ab$ an {\em improper abridged $(a,b)$-word of type II} and refer to a type II abridged $(a,b)$-word $a^nb^j$ with $n>j$ as {\em proper}. In this terminology, if $x\in S$ is represented by an abridged $(a,b)$-word $a^ib^ma^nb^j$ of type III, then $x$ equals the product of an abridged $(a,b)$-word $a^ib^m$ of type I and a proper abridged $(a,b)$-word $a^nb^j$ of type II. 
\bl\label{201} Let $S$ be an orthodox semigroup, $a\in N_S$, $b\in V_S(a)$, and $ab=a^2b^2$. Then
\vspace{0.03in}\\
\indent {\rm(i)} if $x=a^ib^m$ and $x'=a^{i'}b^{m'}$ are abridged words in $a,b$ of type I, then $xx'=a^ib^{m+m'-i'}$, so the product of two abridged $(a,b)$-words of type I equals an abridged $(a,b)$-word of type I;
\vspace{0.03in}\\
\indent{\rm(ii)} if $y=a^nb^j$ and $y'=a^{n'}b^{j'}$ are abridged words in $a,b$ of type II, then $yy'=a^{n'+n-j}b^{j'}$, and hence the product of two abridged words in $a,b$ of type II, at least one of which is proper, equals a proper abridged word of type II;
\vspace{0.03in}\\
\indent{\rm(iii)} if $x\!=\!a^ib^m$ and $y\!=\!a^nb^j$ are abridged words in $a,b$ of types I and II, respectively, then $xy$ equals an abridged word in $a,b$ of type III if $y$ is proper, and $xy\!=\!x$ if $y$ is improper, while $yx\!=\!ab^{m-i+j-n+1}$ if $m-i\!>\!n-j$, and $yx\!=\!a^{n-j+i-m+1}b$ if $m-i\!\leq\! n-j$, so $yx$ equals an abridged word in $a,b$ of type I or II, with the latter being improper when $m-i\!=\!n-j$;
\vspace{0.03in}\\
\indent{\rm(iv)} the product of two abridged words of type III equals an abridged word of type III. 
\el
{\bf Proof.} Statements (i), (ii), and (iii) follow from \cite[Lemma 2.1 and its proof]{key10}. It remains to prove (iv). Suppose that $x$ and $y$ are abridged words in $a,b$ of type III. Then $x=x_1x_2$ and $y=y_1y_2$, where $x_1$ and $y_1$ are abridged $(a,b)$-words of type I, and $x_2$ and $y_2$ are proper abridged $(a,b)$-words of type II. By (iii), $x_2y_1$ is an abridged $(a,b)$-word of type I or II. If $x_2y_1$ is of type I, then according to (i), $x_1(x_2y_1)$ is also of type I, and hence, by (iii), $xy=(x_1x_2y_1)y_2$ is of type III. If $x_2y_1$ is of type II, it follows from (ii) that $(x_2y_1)y_2$ is a proper abridged $(a,b)$-word of type II, and so, by (iii), $xy=x_1(x_2y_1y_2)$ is an abridged $(a,b)$-word of type III. \epr
\vspace{0.1in}\\
\indent Let $S=\langle a,b\rangle$ be a monogenic orthodox semigroup. Suppose that $a\in N_S$ and $ab=a^2b^2$; by \cite[Lemma 1.4]{key10}, this is equivalent to the requirement that $ab=a^2b^2$ and $ba\neq b^2a^2$. Then, according to \cite[Lemma 2.3]{key10}, $S$ is bisimple and combinatorial, and by \cite[Lemma 2.2]{key10}, each $x\in S$ can be written in an abridged form. Recall the following notations introduced in \cite{key10}:
\vspace{0.15in}\\
\indent 1) if $a^n(ab)\not=a^n$ and $(ab)b^n\not=b^n$ for all $n\in\N$, we denote $S$ by $\cO_{(\infty, \infty)}(a,b)$;  
\vspace{0.08in}\\
\indent 2) if $(ab)b^m\not=b^m$ for all $m\in\N$ but $a^k(ab)=a^k$ for some $k\in\N$, then letting $n$ be the smallest of such integers $k$, we denote $S$ by $\cO_{(n, \infty)}(a,b)$;
\vspace{0.08in}\\
\indent 3) if $a^n(ab)\not=a^n$ for all $n\in\N$ but $(ab)b^{l}=b^{l}$ for some $l\in\N$, then with $m$ standing for the smallest of such integers $l$, we denote $S$ by $\cO_{(\infty, m)}(a,b)$; 
\vspace{0.08in}\\
\indent 4) if $a^n(ab)=a^n$ and $(ab)b^m=b^m$ for some $n,m\in\N$, then letting $n$ and $m$ be the smallest integers with these properties, we denote $S$ by $\cO_{(n, m)}(a,b)$ (of course, $\cO_{(1, 1)}(a,b)$ is just another notation for the bicyclic semigroup $\cB(a,b)$). 
\vspace{0.15in}\\
\indent From \cite[Proposition 2.4]{key10}, it follows that if $S=\langle a,b\rangle$ is a monogenic orthodox semigroup such that $a\in N_S$ and $ab=a^2b^2$, then $S=\cO_{(\nu,\,\mu)}(a,b)$ for some $(\nu, \mu)\in\N^{\ast}\times\N^{\ast}$ (as usual, we extend the natural strict linear order $<$ on $\N$ to the set $\N^{\ast}=\N\cup\{\infty\}$ by letting $n<\infty$ for all $n\in\N$). Note that, according to \cite[Remark 4]{key10}, if $S=\cO_{(\nu,\,\mu)}(a,b)$, then $S^{\rm opp}=\cO_{(\mu,\,\nu)}(b,a)$.
 
\br\label{202}\cite[Theorem 2.9]{key10} Let $S=\langle a, b\rangle$ be a monogenic orthodox semigroup where $a$ (and hence $b$) is a nongroup element of $S$. Then $S$ is bisimple if and only if $S$ (or the dual of $S$) coincides with one of the semigroups $\cO_{(\infty, \infty)}(a,b)$, $\cO_{(n, \infty)}(a,b)$, $\cO_{(\infty, m)}(a,b)$, or $\cO_{(n, m)}(a,b)$ for some $m,n\in\N$.
\er

By \cite[Lemmas 2.5, 2.6]{key10}, the eggbox picture of $\cO_{(\infty, \infty)}(a,b)$ is shown in \cite[Fig. 1]{key10}, and those of $\cO_{(n, \infty)}(a,b)$ and $\cO_{(n, m)}(a,b)$ for $m,n\in\N$ are given in \cite[Fig. 2 and Fig. 3]{key10}. For the reader's convenience, we reproduce the eggbox pictures of $\cO_{(\infty, \infty)}(a,b)$, $\cO_{(n, \infty)}(a,b)$, $\cO_{(\infty, m)}(a,b)$, and $\cO_{(n, m)}(a,b)$ for $m, n\in\N$ in Fig. 1.

\bl\label{203} Let $S=\langle a,b\rangle$ be a monogenic orthodox semigroup such that $ab=a^2b^2$ and $a\in N_S$. Then $x\in S$ can be represented by an abridged $(a,b)$-word of type I {\rm[}type II, type III\,{\rm]} if and only if $x\in L_{ab}\setminus\{ab\}$ {\rm[}$x\in R_{ab}$, $x\in S\setminus(L_{ab}\cup R_{ab})${\rm]}.
\el
{\bf Proof.} This is immediate from the eggbox pictures shown in Fig. 1.\epr 
\vspace{0.2in}\\
\[
\renewcommand{\arraystretch}{1.15}
\begin{tabular}{c|c|c|c|c|c|c} 
$\vdots\;\vdots\;\vdots$&\!$\vdots$\!&\!$\vdots$\!&\!$\vdots$\!&\!$\vdots$\!&\!$\vdots$\!&$\vdots\;\vdots\;\vdots$\!\\ \hline
$\cdots$&\!$ab^3a^3b$\!&\!$ab^3a^2b$\!&\!$ab^3$\!&\!$ab^3a$\!&\!$ab^3a^2$\!&$\cdots$\!\\ \hline
$\cdots$&\!$ab^2a^3b$\!&\!$ab^2a^2b$\!&\!$ab^2$\!&\!$ab^2a$\!&\!$ab^2a^2$\!&$\cdots$\!\\ \hline
$\cdots$&\!$a^3b$\! & \!$a^2b$\!  & \!$ab$\! &\!$a$\!&\!$a^2$\!&\!$\cdots$\!\\ \hline
$\cdots$&\!$ba^3b$\! & \!$ba^2b$\! & \!$b$\! &\!$ba$\!&\!$ba^2$\!&$\cdots$\!\\ \hline
$\cdots$&\!$b^2a^3b$\!&\!$b^2a^2b$\!&\!$b^2$\!&\!$b^2a$\!&\!$b^2a^2$\!&$\cdots$\!\\ \hline
$\vdots\;\vdots\;\vdots$&\!$\vdots$\!&\!$\vdots$\!&\!$\vdots$\!&\!$\vdots$\!&\!$\vdots$\!&$\vdots\;\vdots\;\vdots$\!\\ 
\end{tabular} 
\hspace{0.5in}
\begin{tabular}{|c|c|c|c|c|c|c} 
\!$\vdots$\!&$\vdots\;\vdots\;\vdots$&\!$\vdots$\!&\!$\vdots$\!&\!$\vdots$\!&\!$\vdots$\!&$\vdots\;\vdots\;\vdots$\\ \hline
\!$ab^3a^nb$\!&$\cdots$&\!$ab^3a^2b$\!&\!$ab^3$\!&\!$ab^3a$\!&\!$ab^3a^2$\!&$\cdots$\\ \hline
\!$ab^2a^nb$\!&$\cdots$&\!$ab^2a^2b$\!&\!$ab^2$\!&\!$ab^2a$\!&\!$ab^2a^2$\!&$\cdots$\\ \hline
\!$a^nb$\!&$\cdots $&\!$a^2b$\!&\!$ab$\!&\!$a$\!&\!$a^2$\!&$\cdots$\\ \hline
\!$ba^nb$\!&$\cdots $&\!$ba^2b$\!&\!$b$\!&\!$ba$\!&\!$ba^2$\!&$\cdots$\\ \hline
\!$b^2a^nb$\!&$\cdots$&\!$b^2a^2b$\!&\!$b^2$\!&\!$b^2a$\!&\!$b^2a^2$\!&$\cdots$\\ \hline
\!$\vdots$\!&$\vdots\;\vdots\;\vdots$&$\vdots$&\!$\vdots$\!&\!$\vdots$\!&\!$\vdots$\!&$\vdots\;\vdots\;\vdots$\\ 
\end{tabular}  
\]
\vspace{0.03in}\\
\hspace*{1in}(a) $\cO_{(\infty, \infty)}(a,b)$\hspace{2.5in}(b) $\cO_{(n, \infty)}(a,b)$\\ 
\[
\renewcommand{\arraystretch}{1.15}
\!\!\begin{tabular}{c|c|c|c|c|c|c} 
\hline
$\cdots$\!\!&\!$ab^ma^3b$\!&\!$ab^ma^2b$\!&\!$ab^m$\!&\!$ab^ma$\!&\!$ab^ma^2$\!&\!$\cdots$\\ \hline
$\vdots\;\vdots\;\vdots$\!&\!$\vdots$\!&$\vdots$&$\vdots$&$\vdots$&\!$\vdots$\!&\!$\vdots\;\vdots\;\vdots$\\ \hline
$\cdots$\!\!&\!$ab^2a^3b$\!&$ab^2a^2b$&$ab^2$&$ab^2a$&\!$ab^2a^2$\!&\!$\cdots$\\ \hline
$\cdots$\!\!&\!$a^3b$\!&\!$a^2b$\!&\!$ab$\!&\!$a$\!&\!$a^2$\!&\!$\cdots$\\ \hline
$\cdots$\!\!&\!$ba^3b$ \!& $ba^2b$ & $b$ & $ba$ & \!$ba^2$\!&\!$\cdots$\\ \hline
$\cdots$\!\!&\!$b^2a^3b$\!&\!$b^2a^2b$\!&\!$b^2$\!&\!$b^2a$\!&\!$b^2a^2$\!&\!\!$\cdots$\!\\ \hline
$\vdots\;\vdots\;\vdots$\!&\!$\vdots$\!&$\vdots$&$\vdots$&$\vdots$&\!$\vdots$\!&\!$\vdots\;\vdots\;\vdots$\\ 
\end{tabular}
\hspace{0.4in}
\begin{tabular}{|c|c|c|c|c|c|c} 
\hline
\!$ab^ma^nb$\!&\!$\cdots$\!\!&\!$ab^ma^2b$\!&\!$ab^m$\!&\!$ab^ma$\!&\!$ab^ma^2$\!&$\cdots$\\ \hline
\!\!$\vdots$\!&\!$\vdots\;\vdots\;\vdots$\!&\!$\vdots$\!&\!$\vdots$\!&\!\!$\vdots$\!\!&\!$\vdots$\!&$\vdots\;\vdots\;\vdots$\\ \hline
\!\!$ab^2a^nb$\!&\!$\cdots$\!\!&\!$ab^2a^2b$\!&\!$ab^2$\!&\!\!$ab^2a$\!\!&\!$ab^2a^2$\!&$\cdots$\\ \hline
\!\!$a^nb$\!&\!$\cdots$\!\!&\!$a^2b$\!&\!$ab$\!&\!\!$a$\!\!&\!$a^2$\!&$\cdots$\\ \hline
\!\!$ba^nb$\!&\!$\cdots$\!\!&\!$ba^2b$\!&\!$b$\!&\!\!$ba$\!\!&\!$ba^2$\!&$\cdots$\\ \hline
\!\!$b^2a^nb$\!&\!$\cdots$\!\!&\!$b^2a^2b$\!&\!$b^2$\!&\!\!$b^2a$\!\!&\!$b^2a^2$\!&$\cdots$\\ \hline
\!\!$\vdots$\!&\!$\vdots\;\vdots\;\vdots$\!&\!$\vdots$\!&\!$\vdots$\!&\!$\vdots$\!&\!$\vdots$\!&$\vdots\;\vdots\;\vdots$\\ 
\end{tabular} 
\]
\vspace{0.03in}\\
\hspace*{1in}(c) $\cO_{(\infty, m)}(a,b)$\hspace{2.5in}(d) $\cO_{(n, m)}(a,b)$
\vspace{0.005in}\\
\begin{center}Figure 1. The eggbox pictures of semigroups $\cO_{(\nu,\,\mu)}(a,b)$ for $\mu, \nu\in\N^{\ast}$.
\end{center}  
\vspace{0.05in}
\bl\label{204} Let $S=\langle a,b\rangle$ be a monogenic orthodox semigroup such that $ab=a^2b^2$ and $ba\neq b^2a^2$. Then $R_{ab}$, $L_{ab}$, $S\setminus R_{ab}$, and $S\setminus L_{ab}$ are subsemigroups of $S$. Moreover, $R_{ab}=\langle a, ab\rangle$ and $L_{ab}=\langle b, ab\rangle$.
\el
{\bf Proof.} From Lemma \ref{203} and part (ii) [parts (i, iii)] of Lemma \ref{201}, it follows that $R_{ab}$ [$L_{ab}$] is a subsemigroup of $S$ such that $R_{ab}=\langle a, ab\rangle$ [$L_{ab}=\langle b, ab\rangle$].  By Lemma \ref{203}, each element of $S\setminus L_{ab}$ can be represented either by a proper abridged $(a,b)$-word of type II or by an abridged $(a,b)$-word of type III. Let $x,y\in S\setminus L_{ab}$. According to Lemma \ref{201}(ii), if both $x$ and $y$ are represented by proper abridged words of type II, then $xy$ equals a proper abridged word of type II, which implies $xy\in S\setminus L_{ab}$. Suppose $x$ equals an abridged word of type III and $y$ a proper abridged word of type II. Then $x=x_1x_2$ where $x_1$ is an abridged word of type I and $x_2$ a proper abridged word of type II. By Lemma \ref{201}(ii), $x_2y$ equals a proper abridged word of type II, and hence, according to Lemma \ref{201}(iii), $xy=x_1(x_2y)$ equals an abridged word of type III, so that  $xy\in S\setminus L_{ab}$. By Lemma \ref{201}(ii, iii), $yx_1$ has an abridged form of type I or of type II, and $yx=(yx_1)x_2$ can be written either as an abridged word of type III or as a proper abridged word of type II, which means that $yx\in S\setminus L_{ab}$. Finally, if both $x$ and $y$ can be represented by abridged words of type III, then $xy\in S\setminus L_{ab}$ by Lemma \ref{201}(iv). We have shown that $S\setminus L_{ab}$ is a subsemigroup of $S$. By a dual argument, $S\setminus R_{ab}$ is a subsemigroup of $S$ as well. \epr 
\vspace{0.05in}\\
\indent From the eggbox picture of $\cO_{(\infty, \infty)}(a,b)$, it is clear that for each $x\in\cO_{(\infty, \infty)}(a,b)$ there is a unique abridged $(a,b)$-word $a^ib^ma^nb^j$ such that $x=a^ib^ma^nb^j$, and we refer to it as the {\em reduced $(a,b)$-word representing} $x$ (or the {\em reduced form} of $x$). Similarly, for any $m,n\in\N$, the eggbox pictures of $\cO_{(n, \infty)}(a,b)$, $\cO_{(\infty, m)}(a,b)$, and $\cO_{(n, m)}(a,b)$ show that for each $x\in\cO_{(n, \infty)}(a,b)$ [$x\in\cO_{(\infty, m)}(a,b)$, $x\in\cO_{(n, m)}(a,b)$] there is a unique abridged $(a,b)$-word $a^ib^ma^{l}b^j$ with $l\leq n$ [$a^ib^ka^nb^j$ with $k\leq m$, $a^ib^ka^{l}b^j$ with $k\leq m$ and $l\leq n$] such that $x=a^ib^ma^{l}b^j$ [$x=a^ib^ka^nb^j$, $x=a^ib^ka^{l}b^j$], and we call it the {\em reduced $(a,b)$-word representing} $x$ (or the {\em reduced form} of $x$). 
\vspace{0.01in}\\
\indent Let $S=\cO_{(\nu,\,\mu)}(a,b)$ for some $\mu,\nu\in\N^{\ast}$. Suppose that an element $x$ of $S$ is represented by a reduced $(a,b)$-word $a^ib^ma^nb^j$ of type III. Then $x_1=a^ib^m$ is the only reduced $(a,b)$-word of type I and $x_2=a^nb^j$ the only reduced $(a,b)$-word of type II such that $x=x_1x_2$, and we will call $x_1$ and $x_2$ the {\em type I} and  {\em type II components} of $x$, respectively.
\bl\label{205} Let $S\!=\!\langle a,b\rangle$ be a monogenic orthodox semigroup such that $a\!\in\!N_S$ and $ab\!=\!a^2b^2$, so that $S=\cO_{(\nu,\,\mu)}(a,b)$ for some $\mu,\nu\in\N^{\ast}$. Then
\vspace{0.03in}\\
\indent{\rm(i)} $S\setminus R_b$ is a subsemigroup of $S$ if and only if $\mu>1$;
\vspace{0.03in}\\
\indent{\rm(ii)} $S\setminus L_{a}$ is a subsemigroup of $S$ if and only if $\nu>1$. 
\el
{\bf Proof.} Since (ii) is the dual of (i), it is sufficient to prove only (i). If $\mu=1$, then $ab^2=(ab)b=b$ and since $a\not\in R_b$ and $b^2\not\in R_b$, it follows that $S\setminus R_b$ is not a subsemigroup of $S$. Now suppose that $\mu>1$. Take any $x,y\in S\setminus R_b$. Let us show that $xy\in S\setminus R_b$. 
\vspace{0.1in}\\
\indent{\bf Case 1:} $x=a^ib^m$ and $y=a^jb^n$ are reduced words in $a,b$ of type I.
\vspace{0.03in}\\
\indent By Lemma \ref{201}(i), $xy=a^ib^{m+n-j}$. By definition of a type I abridged word, $m-i\geq 1$ and $n-j\geq 1$. If $i=0$, then $xy=b^{m+n-j}\neq b$ since $m+n-j\geq 2$, so that $xy\not\in R_b$. Now suppose that $i=1$ (and thus $m\geq 2$). If $m+n-j-1<\mu$, then $ab^{m+n-j}$ is a reduced word of type I whence $xy=ab^{m+n-j}\not\in R_b$. Finally, if $\mu\in\N$ and $m+n-j=\mu+1$, then $\mu\geq 2$ and therefore $xy=ab^{m+n-j}=(ab)b^{\mu}=b^{\mu}\not\in R_b$. 
\vspace{0.08in}\\
\indent{\bf Case 2:} $x$ and $y$ are both reduced words in $a,b$ of type II.
\vspace{0.03in}\\
\indent In this case, by Lemma \ref{201}(ii), $xy$ is also a reduced word in $a,b$ of type II. Thus $xy\not\in R_b$.
\vspace{0.08in}\\
\indent{\bf Case 3:} $x=a^nb^j$ and $y=a^ib^m$ are reduced words in $a,b$ of types II and I, respectively.
\vspace{0.03in}\\
\indent By Lemma \ref{201}(iii), if $m-i\leq n-j$, then $xy=a^{n-j+i-m+1}b\not\in R_b$ since $a^{n-j+i-m+1}b$ is an abridged word in $a,b$ of type II, and if $m-i>n-j$, then $xy=ab^{m-i+j-n+1}\not\in R_b$ because $ab^{m-i+j-n+1}=b$ only if $m-i=n-j+1$ and $\mu=1$, whereas, by assumption, $\mu>1$.
\vspace{0.08in}\\
\indent{\bf Case 4:} $x=a^ib^m$ and $y=a^nb^j$ are reduced words in $a,b$ of types I and II, respectively.
\vspace{0.03in}\\
\indent
If $y$ is improper, then $xy=x\not\in R_b$. Suppose that $y$ is proper. Then $xy$ is a reduced word of type III whose type I component is $x$ and type II component is $y$. Since $x\not\in R_b$, we have $x\neq b$. Thus $xy\neq ba^nb^j$, which shows that $xy\not\in R_b$.
\vspace{0.08in}\\
\indent{\bf Case 5:} $x$ and $y$ are reduced words in $a,b$ of types I and III, respectively.
\vspace{0.03in}\\
\indent Let $y_1$ and $y_2$ be the type I and type II components of $y$, respectively. Then $xy=(x_1y_1)y_2$ where $x_1y_1$ equals an abridged $(a,b)$-word of type I, and $y_2$ is a proper reduced $(a,b)$-word of type II. Thus $xy$ equals a type III abridged $(a,b)$-word. Let $z_1$ be the reduced form of $x_1y_1$. As shown in Case 1, $z_1\neq b$ because $x_1\neq b$. Since $z_1$ and $y_2$ are the type I and type II components of $xy$, respectively, it follows from Case 4 that $xy\not\in R_b$.
\vspace{0.08in}\\
\indent{\bf Case 6:} $x$ and $y$ are reduced words in $a,b$ of types III and I, respectively.
\vspace{0.03in}\\
\indent Let $x_1$ and $x_2$ be the type I and type II components of $x$, respectively. Then $xy=x_1(x_2y)$ and $x_1\neq b$ since $x\not\in R_b$. By Lemma \ref{201}(iii), $x_2y$ equals an abridged $(a,b)$-word of type I or type II, and if it is of type I, then $xy\not\in R_b$ by Case 1. Suppose $x_2y$ equals an abridged $(a,b)$-word of type II. If that word is improper, then $xy\not\in R_b$ because $xy=x_1\in L_b\setminus\{ab, b\}$, and if it is proper, then $xy\not\in R_b$ by Case 4.
\vspace{0.08in}\\
\indent{\bf Case 7:} $x$ and $y$ are reduced words in $a,b$ of types II and III, respectively.
\vspace{0.03in}\\
\indent Let $y_1$ and $y_2$ be the type I and type II components of $y$, respectively. Then $xy=(xy_1)y_2$. By Lemma \ref{201}(iii) again, $xy_1$ equals an abridged $(a,b)$-word of type I or type II. If the latter is true, then $xy=(xy_1)y_2\not\in R_b$ by Case 2, and if the former holds, then $xy=(xy_1)y_2\not\in R_b$ by Case 4.
\vspace{0.08in}\\
\indent{\bf Case 8:} $x$ and $y$ are reduced words in $a,b$ of types III and II, respectively.
\vspace{0.03in}\\
\indent Let $x_1$ and $x_2$ be the type I and type II components of $x$, respectively. Since $xy=x_1(x_2y)$ and, by Lemma \ref{201}(ii), $x_2y$ equals a proper abridged $(a,b)$-word of type II, we have $xy\not\in R_b$ by Case 4.
\vspace{0.08in}\\
\indent{\bf Case 9:} $x$ and $y$ are both reduced words in $a,b$ of type III.
\vspace{0.03in}\\
\indent Denote the type I and type II components of $x$ by $x_1$ and $x_2$, and those of $y$ by $y_1$ and $y_2$, respectively. Then $xy=x_1(x_2y_1)y_2$. According to Lemma \ref{201}(iii), $x_2y_1$ equals an abridged $(a,b)$-word of type I or of type II. In the former case, $xy_1=x_1(x_2y_1)$ equals an abridged $(a,b)$-word of type I and hence $xy=(xy_1)y_2\not\in R_b$ by Case 4, while in the latter case, by Lemma \ref{201}(ii), $x_2y=(x_2y_1)y_2$ equals a proper abridged $(a,b)$-word of type II, and therefore $xy=x_1(x_2y)\not\in R_b$, again by Case 4.  \epr
\vspace{0.1in}\\
\indent The diagrams of the bands of idempotents of $\cO_{(\infty, \infty)}(a,b)$, $\cO_{(n, \infty)}(a,b)$, and $\cO_{(n, m)}(a,b)$ (with $m\geq n$) are shown in parts (a), (b), and (c), respectively, of \cite[Fig. 4]{key10}. The diagrams of $E_{\cO_{(n, \infty)}(a,b)}$ and $E_{\cO_{(n, m)}(a,b)}$ are drawn in \cite{key10} under the assumptions that $n>1$ and $m>n>1$, respectively, with a remark that modifications for $n=1$ and for $m=n>1$ or $m>n=1$ are obvious. For the reader's convenience, in Fig. 2 we exhibit the diagrams of $E_{\cO_{(\infty,\mu)}(a,b)}$ for $\mu\in\N^{\ast}$, and of $E_{\cO_{(n,m)}(a,b)}$ for $m,n\in\N$ satisfying $n\geq m>1$ or $n>m=1$; as in \cite{key10}, the bold line segments represent the covering relation of the natural order and the thin line segments indicate the $\cR$- and $\cL$-relations on each of these bands. 
\vspace{0.2in}\\
\begin{picture}(520,540)

\put(-5,450){\line(5,2){44}}
\put(-5,396){\line(5,2){44}}
\put(75,450){\line(5,2){44}}
\put(75,396){\line(5,2){44}}

\put(40,504){\line(0,-1){120}}
\put(39.8,504){\line(0,-1){120}}
\put(39.5,504){\line(0,-1){120}}
\put(39.2,504){\line(0,-1){120}}
\put(40.2,504){\line(0,-1){120}}
\put(40.5,504){\line(0,-1){120}}
\put(40.8,504){\line(0,-1){120}}

\put(-5,450){\line(0,-1){85}}
\put(-4.8,450){\line(0,-1){85}}
\put(-4.5,450){\line(0,-1){85}}
\put(-4.2,450){\line(0,-1){85}}
\put(-5.2,450){\line(0,-1){85}}
\put(-5.5,450){\line(0,-1){85}}
\put(-5.8,450){\line(0,-1){85}}

\put(75,450){\line(0,-1){85}}
\put(74.8,450){\line(0,-1){85}}
\put(74.5,450){\line(0,-1){85}}
\put(74.2,450){\line(0,-1){85}}
\put(75.2,450){\line(0,-1){85}}
\put(75.5,450){\line(0,-1){85}}
\put(75.8,450){\line(0,-1){85}}

\put(120,468){\line(0,-1){85}}
\put(119.8,468){\line(0,-1){85}}
\put(119.5,468){\line(0,-1){85}}
\put(119.3,468){\line(0,-1){85}}
\put(120.2,468){\line(0,-1){85}}
\put(120.5,468){\line(0,-1){85}}
\put(120.7,468){\line(0,-1){85}}

\put(-5,450){\line(1,0){80}}
\put(-5,396){\line(1,0){80}}

\put(82, 470){\tiny{$\cR$}}
\put(50, 452){\tiny{$\cR$}}
\put(97, 453){\tiny{$\cL$}}
\put(20, 454){\tiny{$\cL$}}

\put(50,398){\tiny{$\cR$}}
\put(82,416){\tiny{$\cR$}}
\put(101,400){\tiny{$\cL$}}
\put(20,400){\tiny{$\cL$}}

\put(40,468){\line(1,0){80}}
\put(40,414){\line(1,0){80}}

\put(40,504){\circle*{6}}

\put(40,468){\circle*{6}} 
\put(120,468){\circle*{6}} 

\put(40,414){\circle*{6}}
\put(120,414){\circle*{6}}

\put(-5,450){\circle*{6}}
\put(75,450){\circle*{6}}

\put(-5,396){\circle*{6}}
\put(75,396){\circle*{6}}

\put(27,505){\tiny{$ab$}} 

\put(10,468){\tiny{$ab^2a^2b$}}
\put(124,468){\tiny{$ab^2a$}}

\put(12,416){\tiny{$ab^3a^3b$}}
\put(124,414){\tiny{$ab^3a^2$}}

\put(-2,442){\tiny{$ba^2b$}}
\put(80,446){\tiny{$ba$}}
\put(-2,388){\tiny{$b^2a^3b$}}
\put(79,390){\tiny{$b^2a^2$}}

\put(39,365){$\vdots$}
\put(-6,345){$\vdots$}
\put(74,345){$\vdots$}
\put(119,365){$\vdots$}

\put(25,295){(a) \,$E_{\cO_{_{\!(\!\infty,\infty\!)}}\!(a,b)}$}


\put(195,476){\line(5,2){44}}
\put(195,402){\line(5,2){44}}
\put(275,476){\line(5,2){44}}
\put(275,402){\line(5,2){44}}

\put(300,372){\line(0,-1){45}}
\put(300.2,372){\line(0,-1){45}}
\put(300.5,372){\line(0,-1){45}}
\put(300.8,372){\line(0,-1){45}}
\put(299.2,372){\line(0,-1){45}}
\put(299.5,372){\line(0,-1){45}}
\put(299.8,372){\line(0,-1){45}}

\put(219.2,372){\line(0,-1){45}}
\put(219.5,372){\line(0,-1){45}}
\put(219.8,372){\line(0,-1){45}}
\put(220,372){\line(0,-1){45}}
\put(220.8,372){\line(0,-1){45}}
\put(220.5,372){\line(0,-1){45}}
\put(220.2,372){\line(0,-1){45}}

\put(240,524){\line(0,-1){45}}
\put(239.7,524){\line(0,-1){45}}
\put(239.3,524){\line(0,-1){45}}
\put(240.3,524){\line(0,-1){45}}
\put(240.7,524){\line(0,-1){45}}

\put(195,476){\line(0,-1){15}}
\put(194.7,476){\line(0,-1){15}}
\put(194.3,476){\line(0,-1){15}}
\put(195.3,476){\line(0,-1){15}}
\put(195.7,476){\line(0,-1){15}}

\put(275,476){\line(0,-1){15}}
\put(274.7,476){\line(0,-1){15}}
\put(274.3,476){\line(0,-1){15}}
\put(275.3,476){\line(0,-1){15}}
\put(275.7,476){\line(0,-1){15}}

\put(320,494){\line(0,-1){15}}
\put(319.7,494){\line(0,-1){15}}
\put(319.3,494){\line(0,-1){15}}
\put(320.3,494){\line(0,-1){15}}
\put(320.7,494){\line(0,-1){15}}

\put(195,402){\line(0,1){15}}
\put(194.7,402){\line(0,1){15}}
\put(194.3,402){\line(0,1){15}}
\put(195.3,402){\line(0,1){15}}
\put(195.7,402){\line(0,1){15}}

\put(240,420){\line(-2,-5){20}}
\put(240.2,420){\line(-2,-5){20}}
\put(240.5,420){\line(-2,-5){20}}
\put(240.8,420){\line(-2,-5){20}}
\put(239,420){\line(-2,-5){20}}
\put(239.2,420){\line(-2,-5){20}}
\put(239.5,420){\line(-2,-5){20}}
\put(239.8,420){\line(-2,-5){20}}
\put(238.8,420){\line(-2,-5){20}}

\put(195,402){\line(4,-5){24}}
\put(195.2,402){\line(4,-5){24}}
\put(195.5,402){\line(4,-5){24}}
\put(195.8,402){\line(4,-5){24}}
\put(194,402){\line(4,-5){24}}
\put(194.2,402){\line(4,-5){24}}
\put(194.5,402){\line(4,-5){24}}
\put(194.8,402){\line(4,-5){24}}
\put(193.8,402){\line(4,-5){24}}
\put(193.5,402){\line(4,-5){24}}

\put(275,402){\line(0,1){15}}
\put(274.7,402){\line(0,1){15}}
\put(274.3,402){\line(0,1){15}}
\put(275.3,402){\line(0,1){15}}
\put(275.7,402){\line(0,1){15}}

\put(240,420){\line(0,1){15}}
\put(239.7,420){\line(0,1){15}}
\put(239.3,420){\line(0,1){15}}
\put(240.3,420){\line(0,1){15}}
\put(240.7,420){\line(0,1){15}}

\put(320,420){\line(0,1){15}}
\put(319.7,420){\line(0,1){15}}
\put(319.3,420){\line(0,1){15}}
\put(320.3,420){\line(0,1){15}}
\put(320.7,420){\line(0,1){15}}

\put(275,402){\line(4,-5){24}}
\put(274.8,402){\line(4,-5){24}}
\put(274.5,402){\line(4,-5){24}}
\put(274.2,402){\line(4,-5){24}}
\put(275.2,402){\line(4,-5){24}}
\put(275.5,402){\line(4,-5){24}}
\put(275.8,402){\line(4,-5){24}}
\put(276,402){\line(4,-5){24}}

\put(320,419){\line(-2,-5){20}}
\put(319.8,419){\line(-2,-5){20}}
\put(319.5,419){\line(-2,-5){20}}
\put(320.2,419){\line(-2,-5){20}}
\put(320.5,419){\line(-2,-5){20}}
\put(320.8,419){\line(-2,-5){20}}
\put(321,419){\line(-2,-5){20}}
\put(321.2,419){\line(-2,-5){20}}

\put(195,476){\line(1,0){80}}
\put(195,402){\line(1,0){80}}
\put(240,494){\line(1,0){80}}
\put(240,420){\line(1,0){80}}
\put(220,372){\line(1,0){80}}
\put(220,342){\line(1,0){80}}

\put(282,496){\tiny{$\cR$}}
\put(246,478){\tiny{$\cR$}}
\put(287,484.5){\tiny{$\cL$}}
\put(210,485.5){\tiny{$\cL$}}

\put(246,404){\tiny{$\cR$}}
\put(282,422){\tiny{$\cR$}}
\put(287,410.5){\tiny{$\cL$}}
\put(210,411.5){\tiny{$\cL$}}

\put(256,374){\tiny{$\cR$}}
\put(256,344){\tiny{$\cR$}}

\put(240,524){\circle*{6}}

\put(240,494){\circle*{6}} 
\put(320,494){\circle*{6}} 

\put(240,420){\circle*{6}}
\put(320,420){\circle*{6}}

\put(195,476){\circle*{6}}
\put(275,476){\circle*{6}}

\put(195,402){\circle*{6}}
\put(275,402){\circle*{6}}

\put(220,372){\circle*{6}}
\put(220,342){\circle*{6}}

\put(300,372){\circle*{6}}
\put(300,342){\circle*{6}}

\put(227,524){\tiny{$ab$}}
 
\put(210,495){\tiny{$ab^2a^2b$}}
\put(326,495){\tiny{$ab^2a$}}
\put(206,423){\tiny{$ab^ma^mb$}}
\put(324,423){\tiny{$ab^ma^{m-1}$}}
\put(178,370){\tiny{$b^ma^{m+1}b$}}
\put(170,340){\tiny{$b^{m+1}a^{m+2}b$}}

\put(173,476){\tiny{$ba^2b$}}
\put(280,472){\tiny{$ba$}}
\put(153,399){\tiny{$b^{m-1}a^mb$}}
\put(279.4,398.2){\tiny{$b^{m\!-\!1}\!a^{\!m\!-\!1}$}}
\put(306,370){\tiny{$b^ma^m$}}
\put(306,340){\tiny{$b^{m+1}a^{m+1}$}}

\put(239,463){$\vdots$}
\put(239,441){$\vdots$}
\put(194,446){$\vdots$}
\put(194,424){$\vdots$}
\put(274,446){$\vdots$}
\put(274,424){$\vdots$}
\put(319,463){$\vdots$}
\put(319,441){$\vdots$}

\put(219,314){$\vdots$}
\put(299,314){$\vdots$}

\put(195,295){(b) \,$E_{\cO_{_{\!(\!\infty,m\!)}}\!(a,b)}\text{ with }m>1$}


\put(400,524){\line(0,-1){75}}
\put(399.7,524){\line(0,-1){75}}
\put(399.3,524){\line(0,-1){75}}
\put(400.3,524){\line(0,-1){75}}
\put(400.7,524){\line(0,-1){75}}

\put(400,494){\line(1,0){60}}
\put(400,464){\line(1,0){60}}
\put(400,378){\line(1,0){60}}

\put(460,494){\line(0,-1){45}}
\put(459.7,494){\line(0,-1){45}}
\put(459.3,494){\line(0,-1){45}}
\put(460.3,494){\line(0,-1){45}}
\put(460.7,494){\line(0,-1){45}}

\put(400,378){\line(0,1){15}}
\put(399.7,378){\line(0,1){15}}
\put(399.3,378){\line(0,1){15}}
\put(400.3,378){\line(0,1){15}}
\put(400.7,378){\line(0,1){15}}

\put(460,378){\line(0,1){15}}
\put(459.7,378){\line(0,1){15}}
\put(459.3,378){\line(0,1){15}}
\put(460.3,378){\line(0,1){15}}
\put(460.7,378){\line(0,1){15}}

\put(400,378){\line(0,-1){15}}
\put(399.7,378){\line(0,-1){15}}
\put(399.3,378){\line(0,-1){15}}
\put(400.3,378){\line(0,-1){15}}
\put(400.7,378){\line(0,-1){15}}

\put(460,378){\line(0,-1){15}}
\put(459.7,378){\line(0,-1){15}}
\put(459.3,378){\line(0,-1){15}}
\put(460.3,378){\line(0,-1){15}}
\put(460.7,378){\line(0,-1){15}}

\put(400,524){\circle*{6}}
\put(400,494){\circle*{6}} 
\put(400,464){\circle*{6}}
\put(400,378){\circle*{6}}
\put(460,494){\circle*{6}} 
\put(460,464){\circle*{6}}
\put(460,378){\circle*{6}}

\put(386,524){\tiny{$ab$}}
\put(378,494){\tiny{$ba^2b$}}
\put(465,494){\tiny{$ba$}}
\put(374,464){\tiny{$b^2a^3b$}}
\put(465,464){\tiny{$b^2a^2$}}
\put(465,378){\tiny{$b^na^n$}}
\put(363,378){\tiny{$b^na^{n+1}b$}}
 
\put(458.5,408){$\vdots$}
\put(458.5,426){$\vdots$}
\put(458.5,345){$\vdots$}
\put(398.5,408){$\vdots$}
\put(398.5,426){$\vdots$}
\put(398.5,345){$\vdots$}

\put(428,381){\tiny{$\cR$}}
\put(428,467){\tiny{$\cR$}}
\put(428,497){\tiny{$\cR$}}

\put(395,295){(c) \,$E_{\cO_{_{\!(\!\infty,1\!)}}\!(a,b)}$}

\put(25,206){\line(5,2){44}}
\put(25,132){\line(5,2){44}}
\put(105,206){\line(5,2){44}}
\put(105,132){\line(5,2){44}}

\put(130,102){\line(0,-1){15}}
\put(130.2,102){\line(0,-1){15}}
\put(130.5,102){\line(0,-1){15}}
\put(130.8,102){\line(0,-1){15}}
\put(129.2,102){\line(0,-1){15}}
\put(129.5,102){\line(0,-1){15}}
\put(129.8,102){\line(0,-1){15}}

\put(49.2,102){\line(0,-1){15}}
\put(49.5,102){\line(0,-1){15}}
\put(49.8,102){\line(0,-1){15}}
\put(50,102){\line(0,-1){15}}
\put(50.8,102){\line(0,-1){15}}
\put(50.5,102){\line(0,-1){15}}
\put(50.2,102){\line(0,-1){15}}

\put(130,54){\line(0,1){15}}
\put(130.2,54){\line(0,1){15}}
\put(130.5,54){\line(0,1){15}}
\put(130.8,54){\line(0,1){15}}
\put(129.2,54){\line(0,1){15}}
\put(129.5,54){\line(0,1){15}}
\put(129.8,54){\line(0,1){15}}

\put(49.2,54){\line(0,1){15}}
\put(49.5,54){\line(0,1){15}}
\put(49.8,54){\line(0,1){15}}
\put(50,54){\line(0,1){15}}
\put(50.8,54){\line(0,1){15}}
\put(50.5,54){\line(0,1){15}}
\put(50.2,54){\line(0,1){15}}

\put(130,54){\line(-4,-3){40}}
\put(130.2,54){\line(-4,-3){40}}
\put(130.5,54){\line(-4,-3){40}}
\put(130.8,54){\line(-4,-3){40}}
\put(131,54){\line(-4,-3){40}}
\put(131.3,54){\line(-4,-3){40}}
\put(129,54){\line(-4,-3){40}}
\put(129.2,54){\line(-4,-3){40}}
\put(129.5,54){\line(-4,-3){40}}
\put(129.8,54){\line(-4,-3){40}}
\put(128.8,54){\line(-4,-3){40}}

\put(48.5,54){\line(4,-3){40}}
\put(48.8,54){\line(4,-3){40}}
\put(49,54){\line(4,-3){40}}
\put(49.2,54){\line(4,-3){40}}
\put(49.5,54){\line(4,-3){40}}
\put(49.8,54){\line(4,-3){40}}
\put(50,54){\line(4,-3){40}}
\put(51,54){\line(4,-3){40}}
\put(50.8,54){\line(4,-3){40}}
\put(50.5,54){\line(4,-3){40}}
\put(50.2,54){\line(4,-3){40}}

\put(70,254){\line(0,-1){45}}
\put(69.7,254){\line(0,-1){45}}
\put(69.3,254){\line(0,-1){45}}
\put(70.3,254){\line(0,-1){45}}
\put(70.7,254){\line(0,-1){45}}

\put(90,24){\line(0,-1){40}}
\put(89.7,24){\line(0,-1){40}}
\put(89.3,24){\line(0,-1){40}}
\put(90.3,24){\line(0,-1){40}}
\put(90.7,24){\line(0,-1){40}}

\put(25,206){\line(0,-1){15}}
\put(24.7,206){\line(0,-1){15}}
\put(24.3,206){\line(0,-1){15}}
\put(25.3,206){\line(0,-1){15}}
\put(25.7,206){\line(0,-1){15}}

\put(105,206){\line(0,-1){15}}
\put(104.7,206){\line(0,-1){15}}
\put(104.3,206){\line(0,-1){15}}
\put(105.3,206){\line(0,-1){15}}
\put(105.7,206){\line(0,-1){15}}

\put(150,224){\line(0,-1){15}}
\put(149.7,224){\line(0,-1){15}}
\put(149.3,224){\line(0,-1){15}}
\put(150.3,224){\line(0,-1){15}}
\put(150.7,224){\line(0,-1){15}}

\put(25,132){\line(0,1){15}}
\put(24.7,132){\line(0,1){15}}
\put(24.3,132){\line(0,1){15}}
\put(25.3,132){\line(0,1){15}}
\put(25.7,132){\line(0,1){15}}

\put(70,150){\line(-2,-5){20}}
\put(70.2,150){\line(-2,-5){20}}
\put(70.5,150){\line(-2,-5){20}}
\put(70.8,150){\line(-2,-5){20}}
\put(69,150){\line(-2,-5){20}}
\put(69.2,150){\line(-2,-5){20}}
\put(69.5,150){\line(-2,-5){20}}
\put(69.8,150){\line(-2,-5){20}}
\put(68.8,150){\line(-2,-5){20}}

\put(25,132){\line(4,-5){24}}
\put(25.2,132){\line(4,-5){24}}
\put(25.5,132){\line(4,-5){24}}
\put(25.8,132){\line(4,-5){24}}
\put(24,132){\line(4,-5){24}}
\put(24.2,132){\line(4,-5){24}}
\put(24.5,132){\line(4,-5){24}}
\put(24.8,132){\line(4,-5){24}}
\put(23.8,132){\line(4,-5){24}}
\put(23.5,132){\line(4,-5){24}}

\put(105,132){\line(0,1){15}}
\put(104.7,132){\line(0,1){15}}
\put(104.3,132){\line(0,1){15}}
\put(105.3,132){\line(0,1){15}}
\put(105.7,132){\line(0,1){15}}

\put(70,150){\line(0,1){15}}
\put(69.7,150){\line(0,1){15}}
\put(69.3,150){\line(0,1){15}}
\put(70.3,150){\line(0,1){15}}
\put(70.7,150){\line(0,1){15}}

\put(150,150){\line(0,1){15}}
\put(149.7,150){\line(0,1){15}}
\put(149.3,150){\line(0,1){15}}
\put(150.3,150){\line(0,1){15}}
\put(150.7,150){\line(0,1){15}}

\put(105,132){\line(4,-5){24}}
\put(104.8,132){\line(4,-5){24}}
\put(104.5,132){\line(4,-5){24}}
\put(104.2,132){\line(4,-5){24}}
\put(105.2,132){\line(4,-5){24}}
\put(105.5,132){\line(4,-5){24}}
\put(105.8,132){\line(4,-5){24}}
\put(106,132){\line(4,-5){24}}

\put(150,149){\line(-2,-5){20}}
\put(149.8,149){\line(-2,-5){20}}
\put(149.5,149){\line(-2,-5){20}}
\put(150.2,149){\line(-2,-5){20}}
\put(150.5,149){\line(-2,-5){20}}
\put(150.8,149){\line(-2,-5){20}}
\put(151,149){\line(-2,-5){20}}
\put(151.2,149){\line(-2,-5){20}}

\put(25,206){\line(1,0){80}}
\put(25,132){\line(1,0){80}}
\put(70,224){\line(1,0){80}}
\put(70,150){\line(1,0){80}}
\put(50,102){\line(1,0){80}}
\put(50,54){\line(1,0){80}}

\put(108,226){\tiny{$\cR$}}
\put(76,208){\tiny{$\cR$}}
\put(117,214.5){\tiny{$\cL$}}
\put(40,215.5){\tiny{$\cL$}}

\put(76,134){\tiny{$\cR$}}
\put(112,152){\tiny{$\cR$}}
\put(117,140.5){\tiny{$\cL$}}
\put(40,141.5){\tiny{$\cL$}}

\put(86,105){\tiny{$\cR$}}
\put(86,57){\tiny{$\cR$}}

\put(70,254){\circle*{6}}

\put(70,224){\circle*{6}} 
\put(150,224){\circle*{6}} 

\put(70,150){\circle*{6}}
\put(150,150){\circle*{6}}

\put(25,206){\circle*{6}}
\put(105,206){\circle*{6}}

\put(25,132){\circle*{6}}
\put(105,132){\circle*{6}}

\put(50,102){\circle*{6}}
\put(50,54){\circle*{6}}

\put(130,102){\circle*{6}}
\put(130,54){\circle*{6}}

\put(90,24){\circle*{6}}
\put(90,-6){\circle*{6}}

\put(57,254){\tiny{$ab$}}
 
\put(40,225){\tiny{$ab^2a^2b$}}
\put(156,225){\tiny{$ab^2a$}}
\put(36,153){\tiny{$ab^ma^mb$}}
\put(154,153){\tiny{$ab^ma^{m-1}$}}
\put(10,96){\tiny{$b^ma^{m+1}b$}}
\put(15,52){\tiny{$b^{n-\!1}\!a^nb$}}

\put(3,206){\tiny{$ba^2b$}}
\put(110,202){\tiny{$ba$}}
\put(-7,122){\tiny{$b^{m-\!1}\!a^mb$}}
\put(109.4,128.2){\tiny{$b^{m\!-\!1}\!a^{\!m\!-\!1}$}}
\put(134,97){\tiny{$b^ma^m$}}
\put(135,52){\tiny{$b^{n-1}a^{n-1}$}}
\put(96,20){\tiny{$b^na^n$}}
\put(96,-8){\tiny{$b^{n+1}a^{n+1}$}}

\put(69,193){$\vdots$}
\put(69,171){$\vdots$}
\put(24,176){$\vdots$}
\put(24,154){$\vdots$}
\put(104,176){$\vdots$}
\put(104,154){$\vdots$}
\put(149,193){$\vdots$}
\put(149,171){$\vdots$}

\put(49,74){$\vdots$}
\put(129,74){$\vdots$}
\put(89,-28){$\vdots$}

\put(15,-45){(d) $E_{\cO_{_{\!(\!n,m\!)}}\!(a,b)}\text{ with }n>m>1$}


\put(210,206){\line(5,2){44}}
\put(210,132){\line(5,2){44}}
\put(290,206){\line(5,2){44}}
\put(290,132){\line(5,2){44}}

\put(255,254){\line(0,-1){45}}
\put(254.7,254){\line(0,-1){45}}
\put(254.3,254){\line(0,-1){45}}
\put(255.3,254){\line(0,-1){45}}
\put(255.7,254){\line(0,-1){45}}

\put(270,72){\line(0,-1){70}}
\put(269.7,72){\line(0,-1){70}}
\put(269.3,72){\line(0,-1){70}}
\put(270.3,72){\line(0,-1){70}}
\put(270.7,72){\line(0,-1){70}}

\put(210,206){\line(0,-1){15}}
\put(209.7,206){\line(0,-1){15}}
\put(209.3,206){\line(0,-1){15}}
\put(210.3,206){\line(0,-1){15}}
\put(210.7,206){\line(0,-1){15}}

\put(290,206){\line(0,-1){15}}
\put(289.7,206){\line(0,-1){15}}
\put(289.3,206){\line(0,-1){15}}
\put(290.3,206){\line(0,-1){15}}
\put(290.7,206){\line(0,-1){15}}

\put(335,224){\line(0,-1){15}}
\put(334.7,224){\line(0,-1){15}}
\put(334.3,224){\line(0,-1){15}}
\put(335.3,224){\line(0,-1){15}}
\put(335.7,224){\line(0,-1){15}}

\put(210,132){\line(0,1){15}}
\put(209.8,132){\line(0,1){15}}
\put(209.5,132){\line(0,1){15}}
\put(209.2,132){\line(0,1){15}}
\put(210.2,132){\line(0,1){15}}
\put(210.5,132){\line(0,1){15}}
\put(210.8,132){\line(0,1){15}}

\put(255,150){\line(1,-5){15.5}}
\put(255.2,150){\line(1,-5){15.5}}
\put(255.5,150){\line(1,-5){15.5}}
\put(254,150){\line(1,-5){15.5}}
\put(254.2,150){\line(1,-5){15.5}}
\put(254.5,150){\line(1,-5){15.5}}
\put(254.8,150){\line(1,-5){15.5}}
\put(253.8,150){\line(1,-5){15.5}}
\put(253.5,150){\line(1,-5){15.5}}

\put(210,132){\line(1,-1){60}}
\put(210.2,132){\line(1,-1){60}}
\put(210.5,132){\line(1,-1){60}}
\put(210.8,132){\line(1,-1){60}}
\put(209,132){\line(1,-1){60}}
\put(209.2,132){\line(1,-1){60}}
\put(209.5,132){\line(1,-1){60}}
\put(209.8,132){\line(1,-1){60}}
\put(208.8,132){\line(1,-1){60}}
\put(208.5,132){\line(1,-1){60}}

\put(290,132){\line(0,1){15}}
\put(289.8,132){\line(0,1){15}}
\put(289.5,132){\line(0,1){15}}
\put(289.2,132){\line(0,1){15}}
\put(290.2,132){\line(0,1){15}}
\put(290.5,132){\line(0,1){15}}
\put(290.8,132){\line(0,1){15}}

\put(255,150){\line(0,1){15}}
\put(254.8,150){\line(0,1){15}}
\put(254.5,150){\line(0,1){15}}
\put(254.2,150){\line(0,1){15}}
\put(255.2,150){\line(0,1){15}}
\put(255.5,150){\line(0,1){15}}
\put(255.8,150){\line(0,1){15}}

\put(335,150){\line(0,1){15}}
\put(334.8,150){\line(0,1){15}}
\put(334.5,150){\line(0,1){15}}
\put(334.2,150){\line(0,1){15}}
\put(335.2,150){\line(0,1){15}}
\put(335.5,150){\line(0,1){15}}
\put(335.8,150){\line(0,1){15}}

\put(290,132){\line(-1,-3){20}}
\put(289.8,132){\line(-1,-3){20}}
\put(289.5,132){\line(-1,-3){20}}
\put(289.2,132){\line(-1,-3){20}}
\put(290.2,132){\line(-1,-3){20}}
\put(290.5,132){\line(-1,-3){20}}
\put(290.8,132){\line(-1,-3){20}}
\put(291,132){\line(-1,-3){20}}

\put(334.8,149){\line(-5,-6){64}}
\put(335,149){\line(-5,-6){64}}
\put(335.2,149){\line(-5,-6){64}}
\put(335.5,149){\line(-5,-6){64}}
\put(335.8,149){\line(-5,-6){64}}
\put(336,149){\line(-5,-6){64}}
\put(336.2,149){\line(-5,-6){64}}
\put(336.5,149){\line(-5,-6){64}}
\put(336.8,149){\line(-5,-6){64}}
\put(337,149){\line(-5,-6){64}}

\put(210,206){\line(1,0){80}}
\put(210,132){\line(1,0){80}}

\put(294,226){\tiny{$\cR$}}
\put(261,208){\tiny{$\cR$}}
\put(306,215.6){\tiny{$\cL$}}
\put(225,215.6){\tiny{$\cL$}}

\put(261,134){\tiny{$\cR$}}
\put(294,152){\tiny{$\cR$}}
\put(305,141.5){\tiny{$\cL$}}
\put(225,141.5){\tiny{$\cL$}}

\put(255,224){\line(1,0){80}}
\put(255,150){\line(1,0){80}}

\put(255,254){\circle*{6}}

\put(255,224){\circle*{6}} 
\put(335,224){\circle*{6}} 

\put(255,150){\circle*{6}}
\put(335,150){\circle*{6}}

\put(210,206){\circle*{6}}
\put(290,206){\circle*{6}}

\put(210,132){\circle*{6}}
\put(289.3,132){\circle*{6}}

\put(270,72){\circle*{6}}
\put(270,42){\circle*{6}}
\put(270,12){\circle*{6}}

\put(242,254){\tiny{$ab$}}
 
\put(226,227){\tiny{$ab^2a^2b$}}
\put(340,227){\tiny{$ab^2a$}}
\put(224,153){\tiny{$ab^na^nb$}}
\put(338,153){\tiny{$ab^na^{\!n-\!1}$}}

\put(190,200){\tiny{$ba^2b$}}
\put(294,200){\tiny{$ba$}}
\put(180,122){\tiny{$b^{n-\!1}a^nb$}}
\put(292.2,127.3){\tiny{$b^{n\!-\!1}\!\!a^{\!n\!-\!1}$}}
\put(275,68){\tiny{$b^na^n$}}
\put(275,40){\tiny{$b^{n+1}a^{n+1}$}}
\put(275,10){\tiny{$b^{n+2}a^{n+2}$}}

\put(254,193){$\vdots$}
\put(254,171){$\vdots$}
\put(209,176){$\vdots$}
\put(209,154){$\vdots$}
\put(289,176){$\vdots$}
\put(289,154){$\vdots$}
\put(334,193){$\vdots$}
\put(334,171){$\vdots$}

\put(269,-12){$\vdots$}

\put(205,-45){(e) $E_{\cO_{_{\!(\!n,n\!)}}\!(a,b)}\text{ with }n>1$}


\put(400,254){\line(0,-1){75}}
\put(399.7,254){\line(0,-1){75}}
\put(399.3,254){\line(0,-1){75}}
\put(400.3,254){\line(0,-1){75}}
\put(400.7,254){\line(0,-1){75}}

\put(400,224){\line(1,0){60}}
\put(400,194){\line(1,0){60}}
\put(400,108){\line(1,0){60}}

\put(460,224){\line(0,-1){45}}
\put(459.7,224){\line(0,-1){45}}
\put(459.3,224){\line(0,-1){45}}
\put(460.3,224){\line(0,-1){45}}
\put(460.7,224){\line(0,-1){45}}

\put(430,78){\line(0,-1){70}}
\put(429.7,78){\line(0,-1){70}}
\put(429.3,78){\line(0,-1){70}}
\put(430.3,78){\line(0,-1){70}}
\put(430.7,78){\line(0,-1){70}}

\put(400,108){\line(0,1){15}}
\put(399.7,108){\line(0,1){15}}
\put(399.3,108){\line(0,1){15}}
\put(400.3,108){\line(0,1){15}}
\put(400.7,108){\line(0,1){15}}

\put(460,108){\line(0,1){15}}
\put(459.7,108){\line(0,1){15}}
\put(459.3,108){\line(0,1){15}}
\put(460.3,108){\line(0,1){15}}
\put(460.7,108){\line(0,1){15}}

\put(400,108){\line(1,-1){30}}
\put(399.8,108){\line(1,-1){30}}
\put(399.5,108){\line(1,-1){30}}
\put(399.2,108){\line(1,-1){30}}
\put(400.2,108){\line(1,-1){30}}
\put(400.5,108){\line(1,-1){30}}
\put(400.8,108){\line(1,-1){30}}
\put(401,108){\line(1,-1){30}}
\put(399,108){\line(1,-1){30}}
\put(398.8,108){\line(1,-1){30}}

\put(460,108){\line(-1,-1){30}}
\put(459.8,108){\line(-1,-1){30}}
\put(459.5,108){\line(-1,-1){30}}
\put(459.2,108){\line(-1,-1){30}}
\put(459,108){\line(-1,-1){30}}
\put(458.8,108){\line(-1,-1){30}}
\put(460.2,108){\line(-1,-1){30}}
\put(460.5,108){\line(-1,-1){30}}
\put(460.8,108){\line(-1,-1){30}}
\put(461,108){\line(-1,-1){30}}

\put(400,254){\circle*{6}}
\put(400,224){\circle*{6}} 
\put(400,194){\circle*{6}}
\put(400,108){\circle*{6}}
\put(460,224){\circle*{6}} 
\put(460,194){\circle*{6}}
\put(460,108){\circle*{6}}
\put(430,78){\circle*{6}}
\put(430,48){\circle*{6}}
\put(430,18){\circle*{6}}

\put(386,254){\tiny{$ab$}}
\put(378,224){\tiny{$ba^2b$}}
\put(465,224){\tiny{$ba$}}
\put(374,194){\tiny{$b^2a^3b$}}
\put(465,194){\tiny{$b^2a^2$}}
\put(465,108){\tiny{$b^{n-1}a^{n-1}$}}
\put(363,108){\tiny{$b^{n-1}a^nb$}}
\put(436,76){\tiny{$b^na^n$}}
\put(436,46){\tiny{$b^{n+1}a^{n+1}$}}
\put(436,16){\tiny{$b^{n+2}a^{n+2}$}}
 
\put(458.5,138){$\vdots$}
\put(458.5,156){$\vdots$}
\put(398.5,138){$\vdots$}
\put(398.5,156){$\vdots$}
\put(428.5,-5){$\vdots$}

\put(428,111){\tiny{$\cR$}}
\put(428,197){\tiny{$\cR$}}
\put(428,227){\tiny{$\cR$}}

\put(375,-45){(f) $E_{\cO_{_{\!(\!n,1\!)}}\!(a,b)}\text{ with }n>1$}

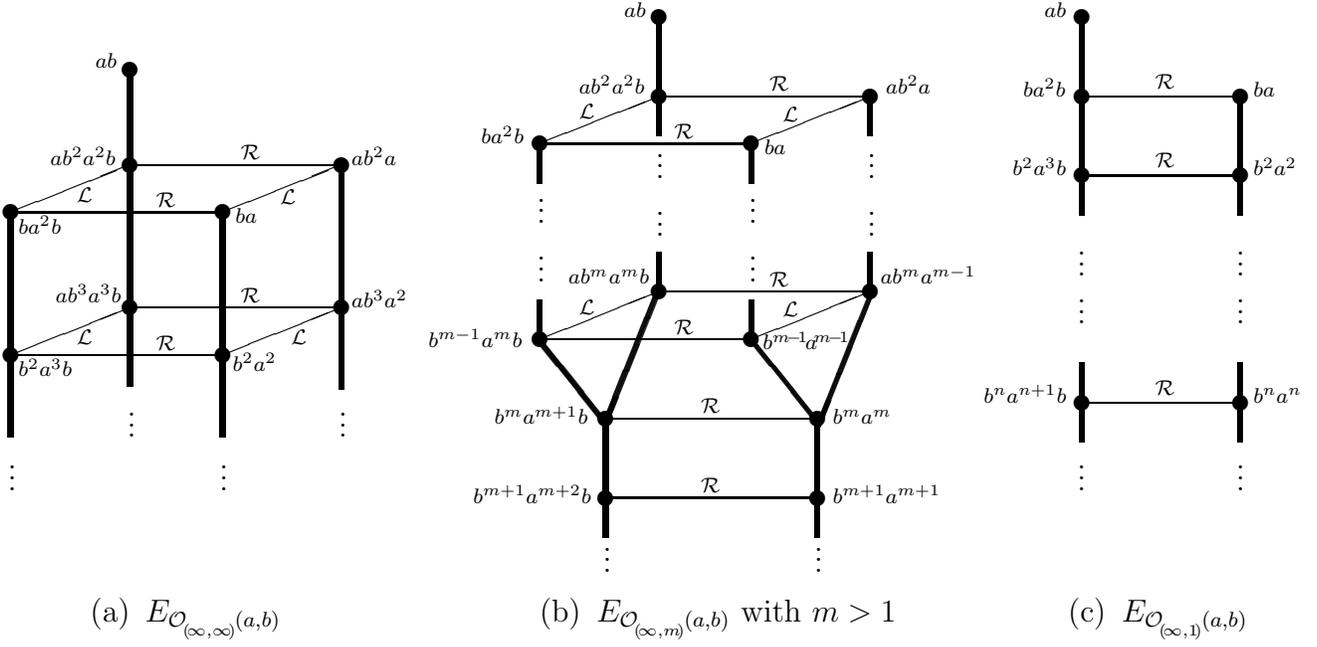
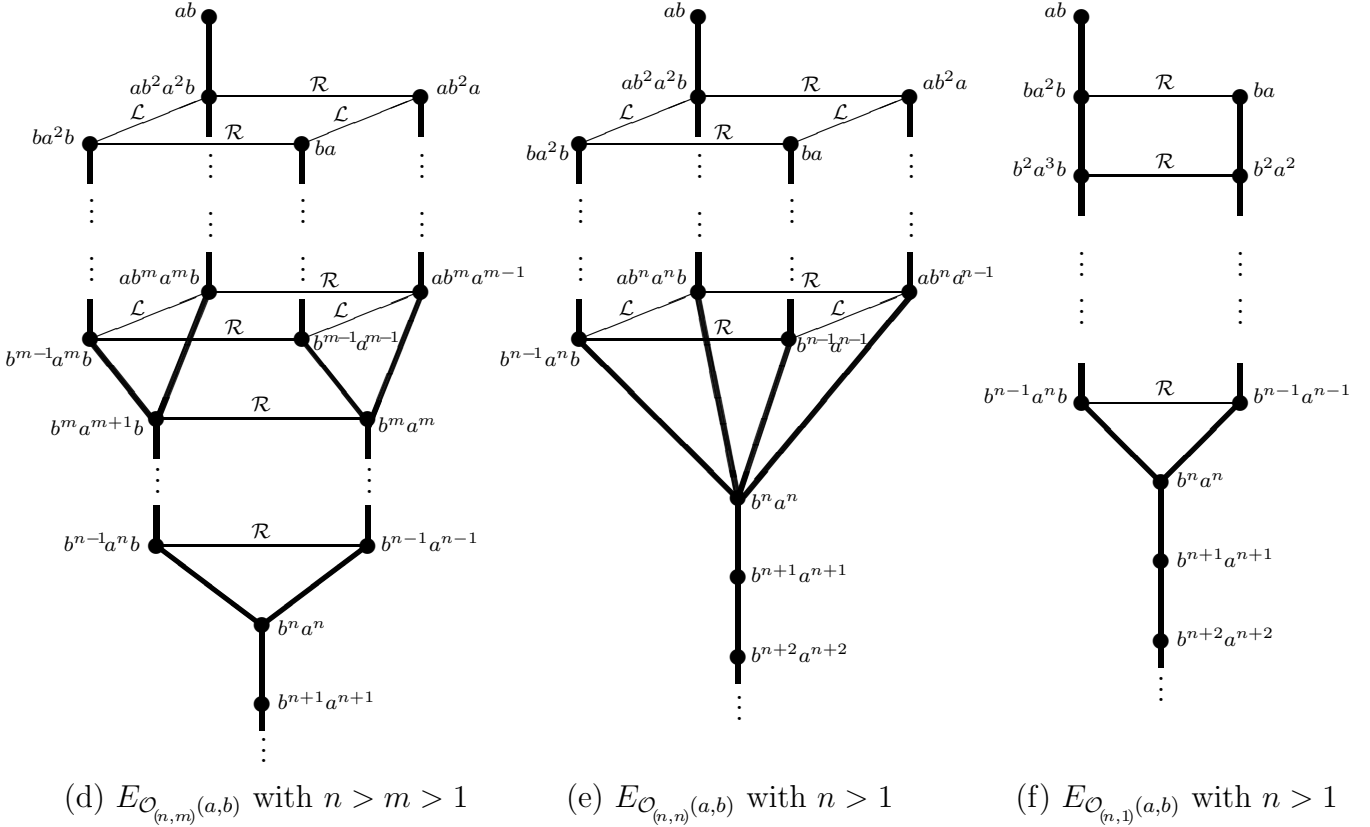
\put(30,-75){$\text{Figure 2. Diagrams of bands of idempotents of }\cO_{(\nu,\,\mu)}(a,b)$ with $\nu\geq\mu>1$ or $\nu>\mu= 1$.}
\end{picture}
\vspace{0.5in}\\
\indent Recall that a band $E$ is {\em uniform} if $eEe\cong fEf$ for all $e,f\in E$. (Since $E$ is equipped with the natural order $\leq$, it is clear that $(eEe,\leq)$ coincides with the principal order ideal $(e\;\!]$ for all $e\in E$.) As shown by Hall \cite[Main Theorem]{key14}, a band $E$ is the band of idempotents of a bisimple orthodox semigroup if and only if $E$ is uniform. In particular, by the easy part of the cited theorem \cite[Result 7]{key14}, the band of idempotents of any bisimple orthodox semigroup is uniform. Thus, by Result \ref{202}, if $E$ is any of the bands shown in Fig. 2, then $E$ is uniform; actually, this is easily seen directly since $eEe$ (for each $e\in E$) is isomorphic to the chain of idempotents of the bicyclic semigroup, so that $eEe\cong fEf$ for all $e,f\in E$. The same remark was made in \cite{key10}, and we would like to take this opportunity to correct an obvious typo on page 321 of \cite{key10} -- on line 4 of the top paragraph on that page (following Fig. 3) it should have been written ``any bisimple orthodox semigroup'' instead of ``any orthodox semigroup''.
\vspace{0.05in}\\
\indent Let $S$ be a regular semigroup. For $x,y\in S$, put $x\leq y$ if and only if $x=ey$ and $x=yf$ for some $e,f\in E_S$; then $\leq$ is an order on $S$ called {\em natural}, and its restriction to $E_S$ coincides with the natural order on $E_S$ (the natural order on $S$ was defined in \cite{key18} by a different condition but it is well known that the definition given here is equivalent to the original one). Recall that $S$ is {\em locally inverse} (or {\em pseudo-inverse}) if $eSe$ is an inverse semigroup for all $e\in E_S$. As shown in \cite{key18}, the following statements are equivalent: (a) $S$ is locally inverse; (b) $\leq$ is compatible with multiplication; (c) for all $e\in E_S$, $(e\;\!]$ is a semilattice. Recall also that a band $E$ is {\em normal} if $efge=egfe$ for all $e,f,g\in E$, and a {\em generalized inverse semigroup} is an orthodox semigroup whose band of idempotents is normal. By \cite[Theorem 1]{key32}, an orthodox semigroup is locally inverse if and only if it is a generalized inverse semigroup. Since for each band $E$ shown in Fig. 2 and for all $e\in E$, the principal order ideal $(e\;\!]$ of $E$ is isomorphic to the chain of idempotents of the bicyclic semigroup, in view of \cite[Theorems 2.9 and 3.1]{key10}, the following fact is immediate.  

\bl\label{206} Every bisimple monogenic orthodox semigroup is locally inverse and hence is a generalized inverse semigroup.
\el

The next lemma, recorded for convenience of reference, can be verified by routine calculation (see also the eggbox pictures in Fig. 1 and the diagrams in Fig. 2 and \cite[Fig. 4]{key10}). 
 
\bl\label{207} Let $S=\cO_{(\nu,\,\mu)}(a,b)$ for some $\mu,\nu\in\N^{\ast}$ such that $\mu>1$ or $\nu>1$. Then
\vspace{0.04in}\\
\indent{\rm(i)} $a^2b\perp ab^2$, $\langle a^2b, ab^2\rangle\!=\!\cB(a^2b, ab^2)$, and the identity element of $\cB(a^2b, ab^2)$ is $ab$, so $ab$ covers $ab^2a^2b$ in the chain of idempotents of $\cB(a^2b, ab^2)$;
\vspace{0.02in}\\
\indent{\rm(ii)} $ab^2a^2\perp ab^3a$, $\langle ab^2a^2, ab^3a\rangle=\cB(ab^2a^2, ab^3a)$, and $ab^2a$ is the identity of $\cB(ab^2a^2, ab^3a)$; 
\vspace{0.02in}\\
\indent{\rm(iii)} $ba^2\perp b^2a$ and $\langle ba^2, b^2a\rangle=\cB(ba^2, b^2a)$ whose identity element is $ba$;
\vspace{0.02in}\\
\indent{\rm(iv)} $ba^3b\perp b^2a^2b$, $\langle ba^3b, b^2a^2b\rangle=\cB(ba^3b, b^2a^2b)$, and the identity of $\cB(ba^3b, b^2a^2b)$ is $ba^2b$;
\vspace{0.02in}\\
\indent{\rm(v)} if $\mu,\nu>1$, then $D^{E_S}_{ab^ka^kb}=\langle ab^ka^kb, b^{k-1}a^{k-1}\rangle=\{ab^ka^kb, ab^ka^{k-1}, b^{k-1}a^{k-1}, b^{k-1}a^kb\}$ is a four-element nonsingular rectangular band for all $k\in\N$ such that $1\leq k-1<{\rm min}\{\mu,\nu\}$;
\vspace{0.02in}\\
\indent{\rm(vi)} if $\mu=1$ {\rm[}$\nu=1${\rm]}, then $\cB(a^2b, ab^2)$ coincides with $\cB(a^2b,b)$ {\rm[}$\cB(a,ab^2)${\rm]}, $\cB(ba^2, b^2a)$ coincides with $\cB(ab^2a^2, ab^3a)$ {\rm[}$\cB(ba^3b, b^2a^2b)${\rm]}, and for all $k\geq 2$ such that $k-1<\nu$ {\rm[}$k-1<\mu${\rm]}, $D^{E_S}_{b^{k-1}a^kb}=\{b^{k-1}a^kb, b^{k-1}a^{k-1}\}$ {\rm[}$D^{E_S}_{ab^ka^{k-1}}=\{ab^ka^{k-1}, b^{k-1}a^{k-1}\}${\rm]} is a two-element right {\rm[}left{\rm]}  singular semigroup.
\el

Let $\mu,\nu\in\N^{\ast}$ be such that $\mu>1$ or $\nu>1$. The following two lemmas can be used for characterizing the semigroup $\cO_{(\nu,\,\mu)}(a,b)$ within the class of bisimple monogenic orthodox semigroups $\langle a, b\rangle$ with nongroup generators by the properties of its subsemigroup lattice. 
 
\bl\label{208} Let $S=\cO_{(\nu,\,\mu)}(a,b)$ for some $\mu,\nu\in\N^{\ast}$ such that $\mu >1$ or $\nu>1$. Then
\vspace{0.05in}\\
\indent{\rm(i)} $S=\cB(a^2b,ab^2)\cup\cB(ba^3b,b^2a^2b)\cup\cB(ba^2,b^2a)\cup\cB(ab^2a^2,ab^3a)\cup\langle a\rangle\cup\langle b\rangle$;
\vspace{0.05in}\\
\indent{\rm(ii)} if $\mu>1$ and $\nu>1$, then the identity elements of the four bicyclic subsemigroups of $S$ listed in {\rm(i)} are pairwise incomparable with respect to the natural order on $E_S$, and if $B$ stands for any of those four bicyclic subsemigroups of $S$ and $C$ is an arbitrary bicyclic subsemigroup of $S$ whose identity element is that of $B$, then $C\subseteq B$;
\vspace{0.05in}\\
\indent{\rm(iii)} if $\nu>\mu>1$, then $\cB(b^{\mu}a^{\mu+2}b, b^{\mu+1}a^{\mu+1}b)$ and $\cB(b^{\mu}a^{\mu+1}, b^{\mu+1}a^{\mu})$ are the largest bicyclic subsemigroups of $\cB(a^2b, ab^2)\cap \cB(ba^3b, b^2a^2b)$ and $\cB(ab^2a^2, ab^3a)\cap \cB(ba^2, b^2a)$, respectively, and, moreover, $\cB(b^{\mu}a^{\mu+2}b, b^{\mu+1}a^{\mu+1}b)\cap\cB(b^{\mu}a^{\mu+1}, b^{\mu+1}a^{\mu})=\emptyset$ when $\nu=\infty$, whereas if $\nu\in\N$, then $\cB(b^{\mu}a^{\mu+2}b, b^{\mu+1}a^{\mu+1}b)\cap\,\cB(b^{\mu}a^{\mu+1}, b^{\mu+1}a^{\mu})\neq\emptyset$ and $\cB(b^{\nu}a^{{\nu}+1}, b^{{\nu}+1}a^{\nu})$ is the largest bicyclic (and proper) subsemigroup of $\cB(b^{\mu}a^{\mu+2}b, b^{\mu+1}a^{\mu+1}b)\cap\cB(b^{\mu}a^{\mu+1}, b^{\mu+1}a^{\mu})$;
\vspace{0.05in}\\
\indent{\rm(iv)} if $\mu>\nu>1$, then $\cB(ab^{\nu+1}a^{\nu+1}, ab^{\nu+2}a^{\nu})$ and $\cB(b^{\nu}a^{\nu+1}, b^{\nu+1}a^{\nu})$ are the largest bicyclic subsemigroups of $\cB(a^2b, ab^2)\cap \cB(ab^2a^2, ab^3a)$ and $ \cB(ba^3b, b^2a^2b)\cap \cB(ba^2, b^2a)$, respectively, and, furthermore, $\cB(ab^{\nu+1}a^{\nu+1}, ab^{\nu+2}a^{\nu})\cap\cB(b^{\nu}a^{\nu+1}, b^{\nu+1}a^{\nu})=\emptyset$ when $\mu=\infty$, and if $\mu\in\N$, then $\cB(ab^{\nu+1}a^{\nu+1}, ab^{\nu+2}a^{\nu})\cap\,\cB(b^{\nu}a^{\nu+1}, b^{\nu+1}a^{\nu})\neq\emptyset$ and $\cB(b^{\mu}a^{{\mu}+1}, b^{{\mu}+1}a^{\mu})$ is the largest bicyclic (and proper) subsemigroup of $\cB(ab^{\nu+1}a^{\nu+1}, ab^{\nu+2}a^{\nu})\cap\cB(b^{\nu}a^{\nu+1}, b^{\nu+1}a^{\nu})$;
\vspace{0.05in}\\
\indent{\rm(v)} if $\nu=\mu\in\N$, then both $\cB(a^2b, ab^2)\cap\,\cB(ba^3b, b^2a^2b)$ and $\cB(ab^2a^2, ab^3a)\cap\,\cB(ba^2, b^2a)$ have the same largest bicyclic subsemigroup, namely, $\cB(b^{\nu}a^{\nu+1}\!, b^{\nu+1}a^{\nu})$. 
\el
{\bf Proof.} Statements (i) and (ii) are immediate from Lemma \ref{207} and the eggbox pictures shown in parts (a), (c), and (d) of Fig. 1.
\vspace{0.06in}\\
\indent{\rm(iii)} Suppose that $\nu>\mu>1$ (in which case, of course, $\mu\in\N$). Assume, first, that $\nu=\infty$. Using part (c) of Fig. 1, we conclude that 
\[\cB(a^2b, ab^2)\cap \cB(ba^3b, b^2a^2b)=\bigcup_{k\geq {\mu}} \!R^{\cB(ba^3b,\, b^2a^2b)}_{b^ka^2b}\text{ and }\cB(ab^2a^2, ab^3a)\cap \cB(ba^2\!, b^2a)=\bigcup_{k\geq {\mu}}\! R^{\cB(ba^2\!,\, b^2a)}_{b^ka}.\]
Take an arbitrary $k\geq {\mu}$. Then $R^{\cB(ba^3b,\, b^2a^2b)}_{b^ka^2b}=\{b^ka^2b,\ldots,b^ka^{\mu}b,b^ka^{{\mu}+1}b,\ldots,b^ka^{k+1}b,\ldots\}$ and the idempotent contained in this $\cR$-class of $\cB(ba^3b,\, b^2a^2b)$ is $b^ka^{k+1}b$. Since 
\[R^{\cB(ba^3b,\, b^2a^2b)}_{b^ka^2b}\cap\cB(b^{\mu}a^{{\mu}+2}b, b^{{\mu}+1}a^{{\mu}+1}b)=\{b^ka^{{\mu}+1}b,\ldots,b^ka^{k+1}b,\ldots\},\] 
it is not difficult to observe that $\cB(b^{\mu}a^{{\mu}+2}b, b^{{\mu}+1}a^{{\mu}+1}b)$ is the largest bicyclic subsemigroup of $\cB(a^2b, ab^2)\cap \cB(ba^3b, b^2a^2b)$. 
Note next that $R^{\cB(ba^2\!,\, b^2a)}_{b^ka}=\{b^ka,\ldots,b^ka^{{\mu}-1},b^ka^{\mu},\ldots,b^ka^k,\ldots\}$ and $b^ka^k$ is the idempotent contained in this $\cR$-class of $\cB(ba^2\!, b^2a)$. Since 
\[R^{\cB(ba^2\!,\, b^2a)}_{b^ka}\cap\cB(b^{\mu}a^{{\mu}+1}\!,\, b^{{\mu}+1}a^{\mu})=\{b^ka^{\mu},\ldots,b^ka^k,\ldots\},\] 
it follows that $\cB(b^{\mu}a^{{\mu}+1}\!, b^{{\mu}+1}a^{\mu})$ is the largest bicyclic subsemigroup of $\cB(ab^2a^2\!, ab^3a)\cap \cB(ba^2\!, b^2a)$. It is also clear that $\cB(b^{\mu}a^{{\mu}+2}b, b^{{\mu}+1}a^{{\mu}+1}b)\cap\cB(b^{\mu}a^{{\mu}+1}\!,\, b^{{\mu}+1}a^{\mu})=\emptyset$ because $\nu=\infty$.

Now assume that $\nu\in\N$. Using part (d) of Fig. 1, we obtain
\[\cB(b^{\mu}a^{{\mu}+2}b, b^{{\mu}+1}a^{{\mu}+1}b)\cap\cB(b^{\mu}a^{{\mu}+1}\!,\, b^{{\mu}+1}a^{\mu})=\bigcup_{k\geq {\nu}}L^{\cB(b^{\mu}a^{{\mu}+1}\!,\, b^{{\mu}+1}a^{\mu})}_{b^{\mu}a^k}.\]
Take an arbitrary $k\geq {\nu}$. Then $L^{\cB(b^{\mu}a^{{\mu}+1}\!,\, b^{{\mu}+1}a^{\mu})}_{b^{\mu}a^k}=\{b^{\mu}a^k,\ldots,b^{{\nu}-1}a^k, b^{\nu}a^k,\ldots, b^ka^k,\ldots\}$ and the idempotent in this $\cL$-class of $\cB(b^{\mu}a^{{\mu}+1}\!, b^{{\mu}+1}a^{\mu})$ is $b^ka^k$. Since 
\[L^{\cB(b^{\mu}a^{{\mu}+1}\!,\, b^{{\mu}+1}a^{\mu})}_{b^{\mu}a^k}\!\cap\cB(b^{\nu}a^{{\nu}+1}\!, b^{{\nu}+1}a^{\nu})\!=\!\{b^{\nu}a^k, \ldots, b^ka^k, \ldots\},\] 
it follows that $\cB(b^{\nu}a^{{\nu}+1}\!, b^{{\nu}+1}a^{\nu})$ is the largest bicyclic (and proper) subsemigroup of the semigroup $\cB(b^{\mu}a^{{\mu}+2}b, b^{{\mu}+1}a^{{\mu}+1}b)\cap\cB(b^{\mu}a^{{\mu}+1}\!, b^{{\mu}+1}a^{\mu})$. 
\vspace{0.05in}\\
\indent The above calculations are based on the eggbox pictures of the semigroups under consideration, constructed using the eggbox pictures shown in parts (c) and (d) of Fig. 1. A visual illustration of statement (iii) is provided by the diagrams in parts (b) and (d) of Fig. 2.
\vspace{0.06in}\\
\indent{\rm(iv)} This statement holds by symmetry with (iii).
\vspace{0.06in}\\
\indent{\rm(v)} This assertion can be checked similarly to the above by examining the eggbox picture in part (d) of Fig. 1 and the diagram in part (e) of Fig. 2. We omit the details.\epr
\vspace{0.06in}\\
\indent If $\nu=1$ or $\mu=1$, then assertion (i) of Lemma \ref{208} is shortened and possible analogs of some other assertions of Lemma \ref{208} can be stated in a simpler form. Due to symmetry, it will be sufficient to consider what happens when $\mu=1$.  
    
\bl\label{209} Let $S=\cO_{(\nu, 1)}(a,b)$ for some $\nu\in\N^{\ast}$ such that $\nu >1$. Then
\vspace{0.05in}\\
\indent{\rm(i)} $S=\cB(a^2b,b)\cup\cB(ba^2,b^2a)\cup\langle a\rangle$;
\vspace{0.03in}\\
\indent{\rm(ii)} the identity elements of $\cB(a^2b,b)$ and $\cB(ba^2,b^2a)$ are incomparable with respect to the natural order on $E_S$, and if $B\in\{\cB(a^2b,b), \cB(ba^2,b^2a)\}$ and $C$ is an arbitrary bicyclic subsemigroup of $S$ whose identity element is that of $B$, then $C\subseteq B$;
\vspace{0.03in}\\
\indent{\rm(iii)} if $\nu\in\N$, then $\cB(b^{\nu}a^{{\nu}+1}, b^{{\nu}+1}a^{\nu})$ is the largest bicyclic (and proper) subsemigroup of $\cB(a^2b, b)\cap\cB(ba^2,b^2a)$.
\el
{\bf Proof.} Statements (i) and (ii) follow easily from Lemma \ref{207}${\rm(i, iv, vi)}$ and the eggbox pictures shown in parts (c) and (d) of Fig. 1 with $m=1$. Note also that since 
\[\cB(a^2b, b)\cap\cB(ba^2,b^2a)=\cB(ba^3b, b^2a^2b)\cap\cB(ba^2,b^2a),\]  
statement (iii) is similar to the situation considered in Lemma \ref{208}(iii) for $\nu\in\N$, except now the condition $\mu>1$ is replaced by $\mu=1$.\epr
\vspace{0.05in}\\
\indent In Section 3, we will have several occasions to use special cases of the following lemma which is easily deduced from the eggbox pictures of the semigroups $\cO_{(\nu,\,\mu)}(a,b)$ shown in Fig. 1. (As usual, the union of an indexed family of sets is empty if the index set is empty.)
\bl\label{210} Let $S=\cO_{(\nu,\,\mu)}(a,b)$ for some $\mu,\nu\in\N^{\ast}$ such that $\mu>1$ or $\nu >1$. Then
\[ \cB(a^2b,ab^2)=S\setminus\left[\left(\bigcup_{1\leq k<\mu}R_{b^k}\right)\cup\left(\bigcup_{1\leq l<\nu}L_{a^l}\right)\right].\]
\el
\bg
\section{Lattice isomorphisms of bisimple monogenic orthodox semigroups}
\sm
In this section we establish the main result of the paper: 

\bt\label{301} Let $S=\langle a,b\rangle$ be an arbitrary bisimple monogenic orthodox semigroup such that $a$ (and hence $b$) has infinite order. Then $S$ is strongly lattice determined.
\et
{\bf Proof}. The proof of this theorem is quite long and it will be convenient to present many of its parts as lemmas. They will be stated and proved at appropriate places within the main proof using previously introduced notations with no additional comments. 
\vspace{0.08in}\\
\indent If $a$ and $b$ are group elements of $S$, then according to \cite[Theorem 3.1]{key10}, $S$ is an $i\times j$ rectangular band of infinite cyclic groups for some $i,j\in\{1, 2\}$, so it is strongly lattice determined by Result \ref{106}. Thus from now on we assume that $a$ (and hence $b$) is a nongroup element of $S$. In view of Result \ref{202}, without loss of generality we will suppose that $S=\cO_{(\nu,\,\mu)}(a,b)$ for some $\mu,\nu\in\N^{\ast}$. If $\mu=\nu=1$, then $S=\cB(a,b)$ and hence it is strongly lattice determined by Result \ref{104}. Therefore {\em through the rest of the proof it will be assumed that $\mu>1$ or $\nu>1$}.  
\vspace{0.08in}\\
\indent Let $\Phi$ be a lattice isomorphism of $S$ onto a semigroup $T$, and let $\varphi$ be the $\Phi$-associated bijection of $S$ onto $T$. By Result \ref{103}, $\varphi|_{E_S}$ is a bijection of $E_S$ onto $E_T$, and each nonidempotent element of $T$ has infinite order. Since $S$ is combinatorial and since the infinite cyclic group is strongly lattice determined, $T$ is combinatorial as well. Denote $x=a\varphi$, $y=b\varphi$, $e=(ab)\varphi$, $u=(a^2b)\varphi$, and $v=(ab^2)\varphi$; clearly, $u\neq x$ if and only if $\nu>1$, and $v\neq y$ if and only if $\mu>1$. Then $T=\langle a, b\rangle\Phi=\langle a\varphi, b\varphi\rangle=\langle x,y\rangle$. By Lemma \ref{207}(i), $\langle a^2b, ab^2\rangle=\cB(a^2b, ab^2)$ and $ab$ is the identity element of $\cB(a^2b, ab^2)$. Since $\langle a^2b, ab^2\rangle\Phi=\langle (a^2b)\varphi, (ab^2)\varphi\rangle=\langle u, v\rangle$, according to Result \ref{104}, $e$ is the identity element of $\langle u,v\rangle$ and $\varphi|_{\cB(a^2b,\,ab^2)}$ is either an isomorphism or an antiisomorphism of $\cB(a^2b,ab^2)$ onto $\langle u,v\rangle$, and if $\varphi|_{\cB(a^2b,\,ab^2)}$ is an isomorphism [antiisomorphism], then $\langle u, v\rangle=\cB(u,v)$ [$\langle u, v\rangle=\cB(v,u)$], which is true if and only if $e=uv$ [$e=vu$]. By Lemma \ref{204}, $L_{ab}$ and $R_{ab}$ are subsemigroups of $S$ and, in fact, $L_{ab}=\langle b, ab\rangle$ and $R_{ab}=\langle a, ab\rangle$. To avoid cumbersome expressions, if $\varphi|_{L_{ab}}$ is an isomorphism [antiisomorphism] of $L_{ab}$ onto $L_{ab}\Phi$, we will often say that it is an isomorphism [antiisomorphism] without mentioning its domain and range, and similarly for $\varphi|_{R_{ab}}$ and $\varphi|_{\cB(a^2b,\,ab^2)}$. Note that if $\nu>1$, then $R_{ab}$ is the semigroup $\langle s, f\rangle$ considered in Proposition \ref{108} with $s=a$ and $f=ab$. It is also obvious that if $\nu=1$, then $R_{ab}=\{ab\}\cup\langle a\rangle$ and $ab$ is an adjoined identity in $R_{ab}$, so if $\varphi|_{R_{ab}}$ is an isomorphism, it is at the same time an antiisomorphism of $R_{ab}$ onto $R_{ab}\Phi$. Likewise, if $\mu=1$, then $L_{ab}=\{ab\}\cup\langle b\rangle$ and if $\varphi|_{L_{ab}}$ is an isomorphism, then it is also an antiisomorphism of $L_{ab}$ onto $L_{ab}\Phi$. {\em The notation and observations of this paragraph will be used, frequently without any reference or explanation, through the rest of the proof.}

\bl\label{302} Suppose that $x=ex$, $y=ye$, and $e=xey$ {\rm[}\;\!$x=xe$, $y=ey$, and $e=yex$\;\!{\rm]}. Then $e=xy$ {\rm[}\;\!$e=yx$\;\!{\rm]} and $T=\langle x, y\rangle$ is a monogenic orthodox semigroup such that $xy=x^2y^2$ and $yx\neq y^2x^2$ {\rm[}\;\!$yx= y^2x^2$ and $xy\neq x^2y^2$\;\!{\rm]}.
\el
{\bf Proof.} Assume that $x=ex$, $y=ye$, and $e=xey$. Then for all $k\in\N$, we have $x^k=ex^k$, $y^k=y^ke$, and $e=x^key^k$, from which it follows that $x^k\cR e\cL y^k$ and so, in particular, $x^k$ and $y^k$ are regular elements of $T$. Therefore all elements of $\langle x\rangle\cup\langle y\rangle$ are regular. By Lemma \ref{208}(i), 
\[S=\cB(a^2b,ab^2)\cup\cB(ba^3b,b^2a^2b)\cup\cB(ba^2,b^2a)\cup\cB(ab^2a^2,ab^3a)\cup\langle a\rangle\cup\langle b\rangle,\] 
and hence, in view of Result \ref{103}, 
\[T=S\varphi=\cB(a^2b,ab^2)\varphi\cup\cB(ba^3b,b^2a^2b)\varphi\cup\cB(ba^2,b^2a)\varphi\cup\cB(ab^2a^2,ab^3a)\varphi\cup\langle x\rangle\cup\langle y\rangle.\]  
Using Result \ref{104}, we conclude that $T$ is a regular semigroup. Moreover, $E_T$ is obviously a subsemigroup of $T$, so $T$ is orthodox. Now $yx(eyxe)yx=y(xey)^2x=yex=yx$ and $(eyxe)yx(eyxe)=ey(xey)^2xe=eyexe=eyxe$, that is, $eyxe\perp yx$. Since $T$ is orthodox and $(eyxe)^2=ey(xey)xe=eyexe=eyxe$, it follows that $yx\in E_T$. As shown above, $x\cR e\cL y$ and since $x=(xe)(yx)$ and $y=(yx)(ey)$, we have $x\cL yx\cR y$. Therefore, by \cite[Theorem 2.17]{key5}, $xy=e$ whence $xyx=x$ and $yxy=y$, that is, $x\perp y$. It follows that $T=\langle x,y\rangle$ is a monogenic orthodox semigroup. Moreover, $x^2y^2=xey=e=xy$ but, according to \cite[Lemma 1.4]{key10}, $yx\neq y^2x^2$ since $x\in N_T$. The alternative statement holds by duality.\epr 
\vspace{0.05in}\\
\indent Applying Proposition \ref{108} to the semigroup $R_{ab}$ when $\nu>1$, and using the fact that $\varphi|_{R_{ab}}$ is both an isomorphism and an antiisomorphism of $R_{ab}$ onto $R_{ab}\Phi$ if $\nu=1$, we obtain 

\bl\label{303} One of the following two statements is true: 
\vspace{0.03in}\\
\indent{\rm(I)} $\varphi|_{R_{ab}}$ is an isomorphism or an antiisomorphism of $R_{ab}$ onto $R_{ab}\Phi$, or 
\vspace{0.03in}\\
\indent{\rm(II)} $xe=ex=u\neq x$,
\vspace{0.03in}\\
with {\rm(II)} being possible only if one of the following holds: either $a^{k+1}b=a^k$ for all $k\geq 2$, which implies $xu=ux=u^2=x^2$ and hence $u^k=x^k$ for all $k\geq 2$, or $a^3b\neq a^2$ but $a^{k+1}b=a^k$ for all $k\geq 3$, in which case $xu=ux=u^2\neq x^2$ and $u^k=x^k$ for all $k\geq 3$.
\el
Although in the case when $\nu>1$, Lemma \ref{303} is just a restatement of Proposition \ref{108} for the semigroup $R_{ab}$, we have omitted from the formulation of that lemma the assertion that $e$ is an adjoined identity in $\langle e,u\rangle$ since it is an immediate consequence of the fact that $e$ is the identity of the bicyclic semigroup $\langle u,v\rangle$. A similar remark can be made about the next lemma obtained by applying the dual of Proposition \ref{108} to the semigroup $L_{ab}$ when $\mu>1$.

\bl\label{304} One of the following two statements holds:
\vspace{0.03in}\\
\indent{\rm(I)} $\varphi|_{L_{ab}}$ is an isomorphism or an antiisomorphism of $L_{ab}$ onto $L_{ab}\Phi$, or 
\vspace{0.03in}\\
\indent{\rm(II)} $ye=ey=v\neq y$,
\vspace{0.03in}\\
with {\rm(II)} being possible only if one of the following is true: either $ab^{k+1}=b^k$ for all $k\geq 2$, which implies $yv=vy=v^2=y^2$ and hence $v^k=y^k$ for all $k\geq 2$, or $ab^3\neq b^2$ but $ab^{k+1}=b^k$ for all $k\geq 3$, in which case $yv=vy=v^2\neq y^2$ and $v^k=y^k$ for all $k\geq 3$.
\el
It must be emphasized that the only facts about $\Phi$ used for obtaining Lemmas \ref{303} and \ref{304} were those asserting that $\Phi|_{\cSb(R_{ab})}$ and $\Phi|_{\cSb(L_{ab})}$ are lattice isomorphisms of $R_{ab}$ onto $R_{ab}\Phi$ and of $L_{ab}$ onto $L_{ab}\Phi$, respectively. However, since $\Phi$ is a lattice isomorphism of $S$ onto $T$, more information is available, and we will eventually be able to show that statement (II) in each of these lemmas is not an actual possibility. 
\vspace{0.05in}\\
\indent Although parts (i) and (ii) of the next lemma (established by routine calculation) are symmetric, for clarity we explicitly state both of them.

\bl\label{305} {\rm(i)} If $\varphi|_{L_{ab}}$ is an isomorphism {\rm[}antiisomorphism\;\!{\rm]} of $L_{ab}$ onto $L_{ab}\Phi$, then  $ye=y$, $ey=v$, $yv=y^2$, and $vy=v^2$ {\rm[}\;\!$ye=v$, $ey=y$, $yv=v^2$, and $vy=y^2$\;\!{\rm]}. 
\vspace{0.03in}\\
\indent{\rm(ii)} If $\varphi|_{R_{ab}}$ is an isomorphism {\rm[}antiisomorphism\;\!{\rm]} of $R_{ab}$ onto $R_{ab}\Phi$, then  $xe=u$, $ex=x$, $xu=u^2$, and $ux=x^2$ {\rm[}\;\!$xe=x$, $ex=u$, $xu=x^2$, and $ux=u^2$\;\!{\rm]}.    
\el

\bl\label{306}  If $\mu>1$ and $\varphi|_{L_{ab}}$ is an isomorphism {\rm[}antiisomorphism\;\!{\rm]} of $L_{ab}$ onto $L_{ab}\Phi$, or if $\nu>1$ and $\varphi|_{R_{ab}}$ is an isomorphism {\rm[}antiisomorphism\;\!{\rm]} of $R_{ab}$ onto $R_{ab}\Phi$, then $e=uv$ {\rm[}\;\!$e=vu$\;\!{\rm]}.
\el
{\bf Proof.} Let $\mu>1$. Suppose $\varphi|_{L_{ab}}$ is an isomorphism of $L_{ab}$ onto $L_{ab}\Phi$. By Lemma \ref{305}(i), $ye=y$ and $yv=y^2$. We know that $e$ is the identity element of $\langle u,v\rangle$ and either $\langle u, v\rangle=\cB(u,v)$ or $\langle u, v\rangle=\cB(v,u)$. Suppose $\langle u, v\rangle=\cB(v,u)$. Then $e=vu$, so that $y=yvu=y^2u\in\langle y^2, u\rangle$ whence $b\in\langle b^2, a^2b\rangle$. By Lemma \ref{205}(i), $S\setminus R_b$ is a subsemigroup of $S$. Since $b^2, a^2b\in S\setminus R_b$, it follows that $b\in\langle b^2, a^2b\rangle\subseteq S\setminus R_b$; a contradiction. Therefore $\langle u, v\rangle=\cB(u,v)$ and so $e=uv$. By duality, if $\varphi|_{L_{ab}}$ is an antiisomorphism of $L_{ab}$ onto $L_{ab}\Phi$, then $e=vu$.

If $\nu>1$ and $\varphi|_{R_{ab}}$ is an isomorphism [antiisomorphism] of $R_{ab}$ onto $R_{ab}\Phi$, then using Lemmas \ref{305}(ii) and \ref{205}(ii), we establish that $e=uv$ [$e=vu$] by a symmetric argument.   \epr

\bl\label{307} Suppose that $\varphi|_{L_{ab}}$ is an isomorphism or an antiisomorphism of $L_{ab}$ onto $L_{ab}\Phi$, and $\varphi|_{R_{ab}}$ is an isomorphism or an antiisomorphism of $R_{ab}$ onto $R_{ab}\Phi$. Then either both $\varphi|_{L_{ab}}$ and $\varphi|_{R_{ab}}$ are isomorphisms or they both are antiisomorphisms.
\el
{\bf Proof.} Assume that $\nu>1$. Suppose that $\varphi|_{R_{ab}}$ is an isomorphism [antiisomorphism] of $R_{ab}$ onto $R_{ab}\Phi$. Then, according to Lemma \ref{306}, $e=uv$ [$e=vu$] whence $\langle u,v\rangle=\cB(u,v)$ and $vu\neq e$ [$\langle u,v\rangle=\cB(v,u)$ and $uv\neq e$]. If $\mu=1$, then it is trivial that $\varphi|_{L_{ab}}$ is an isomorphism [antiisomorphism] of $L_{ab}$ onto $L_{ab}\Phi$. Suppose that $\mu>1$. If $\varphi|_{L_{ab}}$ were an antiisomorphism [isomorphism] of $L_{ab}$ onto $L_{ab}\Phi$, by Lemma \ref{306} again, we would have $e=vu$ [$e=uv$]; a contradiction. Therefore $\varphi|_{L_{ab}}$ is an isomorphism [antiisomorphism] of $L_{ab}$ onto $L_{ab}\Phi$. 

By symmetry, if $\mu>1$ and $\varphi|_{L_{ab}}$ is an isomorphism [antiisomorphism] of $L_{ab}$ onto $L_{ab}\Phi$, then $\varphi|_{R_{ab}}$ is an isomorphism [antiisomorphism] of $R_{ab}$ onto $R_{ab}\Phi$.  \epr
\vspace{0.08in}\\
\indent As noted earlier, we intend to show that in Lemmas \ref{303} and \ref{304} only statement (I) is actually true; in view of duality, we may do this only for Lemma \ref{303}. To achieve this goal, we will prove several auxiliary results. It will be sufficient to show that statement (II) of Lemma \ref{303} does not hold under the assumption that $e=uv$; if $e=vu$, the same assertion is established by a dual argument. If $e=uv$, then according to Lemmas \ref{306} and \ref{304}, either $\varphi|_{L_{ab}}$ is an isomorphism of $L_{ab}$ onto $L_{ab}\Phi$, or statement (II) of Lemma \ref{304} holds. Recall also that if $e=uv$, then $\langle u,v\rangle=\cB(u,v)$ and $\varphi|_{\cB(a^2b,\,ab^2)}$ is an isomorphism of $\cB(a^2b,ab^2)$ onto $\cB(u,v)$. {\em In what follows the observations of this paragraph will be used without mention.}

\bl\label{308} Suppose that statement {\rm(II)} of Lemma {\rm\ref{303}} is true and $e=uv$. Then $e=xy$ and $yx\in E_T\setminus\{e\}$.
\el
{\bf Proof.} By Lemma \ref{303}, $xe=ex=u\neq x$ and either $\nu=2$ (whence $xu=ux=u^2=x^2$ and so $u^k=x^k$ for all $k\geq 2$), or $\nu=3$ (in which case $xu=ux=u^2\neq x^2$ and $u^k=x^k$ for all $k\geq 3$); these facts will be used below with no comments. If $\mu=1$, then $v=y$ and since $ev=ve=v$ and $vu\neq e$, we have $e=uv=xey=xy$ and $(yx)^2=y(xy)x=yex=yx\neq e$. Thus for the rest of the proof it will be assumed that $\mu>1$ (and so $y\neq v$). 
\vspace{0.08in}\\
\indent {\bf Case 1:} $\varphi|_{L_{ab}}$ is an isomorphism of $L_{ab}$ onto $L_{ab}\Phi$.
\vspace{0.06in}\\
\indent By Lemma \ref{305}(i), $ye=y$, $ey=v$, $yv=y^2$, and $vy=v^2$. Hence $e=uv=xey=exy$ which implies $y=ye=y(exy)=(ye)xy=yxy$. Therefore $yx=(yx)^2$ and $xy=(xy)^2$. 
\vspace{0.06in}\\
\indent Suppose that $e\neq xy$. Since $(xy)e=xy$ and $e(xy)=xey=uv=e$, we conclude that $\{e,xy\}$ is a two-element left singular subsemigroup of $T$. Thus $\{ab,(xy)\varphi^{-1}\}$ is either a two-element singular semigroup or a two-element semilattice. However $ab$ is not contained in any two-element singular subsemigroup of $S$. Hence  $\{ab,(xy)\varphi^{-1}\}$ is a two-element semilattice. To determine the possible values of $(xy)\varphi^{-1}$, we will list the elements of the principal order ideal $(ab\;\!]$ of $E_{\cO_{(\nu,\,\mu)}(a,b)}$. From the diagrams in parts (b) and (c) of \cite[Fig. 4]{key10} and those in parts (d) and (e) of Fig. 2, we deduce that if $\nu=2$, then

\setcounter{equation}{8}
\be
(ab\;\!]\!=\!\left\{\!\! 
\begin{array}{l}\{ab, ab^2a^2b\}\cup\{b^k a^k : k\geq 2\} \\ \{ab, ab^2a^2b\}\cup\{ab^{l+1}a^l:2\leq l<\mu\}\cup\{b^k a^k:k\geq \mu\} \\ \{ab, ab^2a^2b\}\cup\{ab^{l+1}a^l:l\geq 2\} \end{array}
\begin{array}{l} \!\!\!\text{ if }\mu=2, \\ \!\!\!\text{ if }2<\mu<\infty,\\ \!\!\!\text{ if }\mu=\infty,
				\end{array}\right.
\label{eq:4th}
\ee 
and if $\nu=3$, then
\be
(ab\;\!]\!=\!\left\{\!\! 
\begin{array}{l}\{ab,\, ab^2a^2b,\, b^2a^3b\;\!\}\cup\{b^k a^k : k\geq 3\} \\ \{ab, ab^2a^2b, ab^3a^3b\}\cup\{b^k a^k : k\geq 3\} \\ \{ab, ab^2a^2b, ab^3a^3b\}\cup\{ab^{l+1}a^l:3\leq l<\mu\}\cup\{b^k a^k:k\geq \mu\} \\ \{ab, ab^2a^2b, ab^3a^3b\}\cup\{ab^{l+1}a^l:l\geq 3\} \end{array}
\begin{array}{l} \!\!\!\text{ if }\mu=2, \\ \!\!\!\text{ if }\mu=3, \\ \!\!\!\text{ if }3<\mu<\infty, \\ \!\!\!\text{ if }\mu=\infty.
				\end{array}\right.
\label{eq:5th}
\ee
Since $\{ab,(xy)\varphi^{-1}\}$ is a two-element semilattice, it follows that $(xy)\varphi^{-1}\in(ab\;\!]\setminus\{ab\}$, so the possible values of $(xy)\varphi^{-1}$ can be found using (\ref{eq:4th}) and (\ref{eq:5th}). In fact, by examining (\ref{eq:4th}) and (\ref{eq:5th}), all possibilities for $(xy)\varphi^{-1}$ can be divided into the following five cases.
\vspace{0.05in}\\
\indent 1(a) $(xy)\varphi^{-1}\!=ab^2a^2b$. 
\vspace{0.02in}\\ 
\indent In this case, $xy=(ab^2a^2b)\varphi=(ab^2)\varphi\cdot(a^2b)\varphi=vu=(ey)(ex)=eyx$ and therefore $y=y(xy)=y(eyx)=y^2x\in\langle y^2, x\rangle$. It follows that $b\in\langle b^2, a\rangle$. Since $b^2,a\in S\setminus R_b$, using Lemma \ref{205}(i), we obtain $\langle b^2, a\rangle\subseteq S\setminus R_b$ whence $b\in S\setminus R_b$; a contradiction.
\vspace{0.05in}\\
\indent 1(b) $(xy)\varphi^{-1}=b^2a^3b$ (when $\mu=2$ and $\nu=3$). 
\vspace{0.02in}\\
\indent  Since $\mu=2$, it is immediate that $b^2=ab^3=(ab^2)^2$. Therefore
\[xy=(b^2a^3b)\varphi=(b^2)\varphi\cdot(a^3b)\varphi=y^2u^2=y^2ex^2=y^2x^2.\] 
It follows that $y=yxy=y^3x^2\in\langle y^3, x^2\rangle$ whence $b\in\langle b^3, a^2\rangle$. Since $b^3,a^2\in S\setminus R_b$, according to Lemma \ref{205}(i), we have $b\in \langle b^3, a^2\rangle\subseteq S\setminus R_b$; a contradiction. 
\vspace{0.05in}\\
\indent 1(c) $(xy)\varphi^{-1}=ab^3a^3b$ (when $\mu\geq \nu=3$). 
\vspace{0.02in}\\ 
\indent Note that $ab^3=(ab^2)^2\in\cB(a^2b,\,ab^2)$ and $a^3b=(a^2b)^2\in\cB(a^2b,\,ab^2)$. Therefore 
\[xy=(ab^3a^3b)\varphi=(ab^3)\varphi\cdot(a^3b)\varphi=v^2u^2=(ey^2)(ex^2)=ey^2x^2.\]
Hence $y=y(ey^2x^2)=y^3x^2\in\langle y^3, x^2\rangle$, which leads to a contradiction exactly as in 1(b).
\vspace{0.05in}\\
\indent 1(d) $(xy)\varphi^{-1}=ab^{l+1}a^l$ for some $l\in\N$ such that $\mu>l\geq\nu$ (when $\mu>\nu$). 
\vspace{0.02in}\\ 
\indent In this case, $a^l=(a^2b)^l\in\cB(a^2b,ab^2)$ and $ab^{l+1}=(ab^2)^l\in\cB(a^2b,ab^2)$. It follows that $xy=(ab^{l+1})\varphi\cdot(a^l)\varphi=v^lx^l=ey^lx^l$ whence $y=y(xy)=y(ey^lx^l)=y^{l+1}x^l\in\langle y^{l+1},x^l\rangle$, which implies $b\in\langle b^{l+1}, a^l\rangle$. Observing that $b^{l+1}, a^l\in S\setminus R_b$ and using Lemma  \ref{205}(i), we conclude that $\langle b^{l+1}, a^l\rangle\subseteq S\setminus R_b$ and therefore $b\in S\setminus R_b$; a contradiction.  
\vspace{0.05in}\\
\indent 1(e) $(xy)\varphi^{-1}=b^ka^k$ for some $k\geq\max\{\nu,\mu\}$.
\vspace{0.02in}\\ 
\indent Here $xy=(b^ka^k)\varphi=(b^k)\varphi\cdot(a^k)\varphi=y^kx^k$ whence $y=yxy=y^{k+1}x^k\in\langle y^{k+1},x^k\rangle$ and therefore $b\in\langle b^{k+1}, a^k\rangle$. Since $b^{k+1}, a^k\notin R_b$, by Lemma \ref{205}(i) we have $b\in\langle b^{k+1}, a^k\rangle\subseteq S\setminus R_b$; again a contradiction.
\vspace{0.06in}\\
\indent We have shown that the assumption that $e\neq xy$ leads to a contradiction in all possible cases. It follows that $e=xy$. As noted earlier, $yx\in E_T$. If $yx=e$, then $y=yxy=ey=v$, which is not true. Therefore $yx\neq e$. 
\vspace{0.06in}\\
\indent {\bf Case 2:} Statement (II) of Lemma \ref{304} holds.
\vspace{0.04in}\\
\indent By Lemma \ref{304}, $ye=ey=v\, (\neq y)$ and either $\mu=2$, which implies $yv=vy=v^2=y^2$  and so $v^k=y^k$ for all $k\geq 2$, or $\mu=3$, in which case $yv=vy=v^2\neq y^2$  and $v^k=y^k$ for all $k\geq 3$ (these observations will be used below without mention). 
\vspace{0.02in}\\
\indent Suppose that $e\neq xy$. By Lemma \ref{210}, 
\be
\cB(a^2b, ab^2)=\left\{\!\! 
\begin{array}{l} S\setminus(R_b\cup L_a)\\ S\setminus(R_b\cup L_a\cup L_{a^2})\\ S\setminus(R_b\cup R_{b^2}\cup L_a)\\ S\setminus(R_b\cup R_{b^2}\cup L_a\cup L_{a^2})\end{array}
\begin{array}{l} \!\!\text{ if }\mu=\nu=2, \\ \!\!\text{ if }\mu=2\text{ and }\nu=3, \\ \!\!\text{ if }\mu=3\text{ and }\nu=2, \\ \!\!\text{ if }\mu=\nu=3.
				\end{array}\right.
\label{eq:6th}
\ee
If $xy\in\cB(u,v)$, then $xy=e(xy)=xey=uv=e$. Since, by assumption, $e\neq xy$, we conclude that $xy\not\in\cB(u,v)$ and therefore $(xy)\varphi^{-1}\not\in\cB(a^2b, ab^2)$. From this and (\ref{eq:6th}), it follows that  
\be
(xy)\varphi^{-1}\in\left\{\!\! 
\begin{array}{l}R_b\cup L_a\\ R_b\cup L_a\cup L_{a^2}\\ R_b\cup R_{b^2}\cup L_a\\ R_b\cup R_{b^2}\cup L_a\cup L_{a^2}\end{array}
\begin{array}{l} \!\!\text{ if }\mu=\nu=2, \\ \!\!\text{ if }\mu=2\text{ and }\nu=3, \\ \!\!\text{ if }\mu=3\text{ and }\nu=2, \\ \!\!\text{ if }\mu=\nu=3.
				\end{array}\right.
\label{eq:7th}
\ee
\indent Suppose that $(xy)\varphi^{-1}\in L_a$.  The eggbox pictures of $\cO_{(\nu,\,\mu)}(a,b)$, where $\nu,\mu\in\{2,3\}$,  show that $L_a=\{ab^2a, a\}\cup\{b^ka:k\in\N\}$ if $\mu=2$, and $L_a=\{ab^3a, ab^2a, a\}\cup\{b^ka:k\in\N\}$ if $\mu=3$, and it is immediate that in both cases, $L_a\cap E_S=\{ab^2a, ba\}$. 

Assume that $xy\in E_T$. Then $(xy)\varphi^{-1}\in\{ab^2a, ba\}$. Since $e(xy)=e=(xy)e$, it follows that $\langle e, xy\rangle=\{e, xy\}$ and hence $|\langle e, xy\rangle|=2$. However $|\langle ab, ba\rangle|=|\{ab, ba, ab^2a, ba^2b, ab^2a^2b\}|=5$ and $|\langle ab, ab^2a\rangle|=|\{ab, ab^2a, ab^2a^2b\}|=3$, so that $|\langle e, xy\rangle|=|\langle ab, ba\rangle|=5$  if $(xy)\varphi^{-1}=ba$, and $|\langle e, xy\rangle|=|\langle ab, ab^2a\rangle|=3$  if $(xy)\varphi^{-1}=ab^2a$. This contradiction shows that $xy\notin E_T$.

Since  $xy\notin E_T$ and $e(xy)=e=(xy)e$, it follows that $o(xy)=\infty$ and $e$ is an adjoined identity in the subsemigroup $\langle e, xy\rangle$ of $T$. Therefore, by \cite[Lemma 33.4]{key30}, $ab$ is either an adjoined identity or an adjoined zero in the subsemigroup $\langle ab, (xy)\varphi^{-1}\rangle$ of $S$. At the same time, it is easily seen that $ab$ is neither a zero nor an identity for any element of $L_a\setminus\{ba, ab^2a\}$, and since $(xy)\varphi^{-1}\in L_a\setminus\{ba, ab^2a\}$, we again have a contradiction. 

We have shown that the assumption that $(xy)\varphi^{-1}\in L_a$ leads to a contradiction in all possible cases. Therefore  $(xy)\varphi^{-1}\notin L_a$, and hence, by symmetry, $(xy)\varphi^{-1}\notin R_b$.
\vspace{0.05in}\\
\indent Now let $(xy)\varphi^{-1}\in L_{a^2}$ when $\nu=3$. Using the eggbox pictures of $\cO_{(3,\mu)}(a,b)$ with $\mu\in\{2,3\}$, we see that if $\mu=2$, then $L_{a^2}=\{ab^2a^2, a^2\}\cup\{b^ka^2:k\in\N\}$ and $L_{a^2}\cap E_S=\{b^2a^2\}$, whereas if $\mu=3$, then $L_{a^2}=\{ab^3a^2, ab^2a^2, a^2\}\cup\{b^ka^2:k\in\N\}$ and $L_{a^2}\cap E_S=\{ab^3a^2, b^2a^2\}$. 

Suppose that $xy\in E_T$. As above, using the fact that $e(xy)=e=(xy)e$, we conclude that $\langle e, xy\rangle=\{e, xy\}$ and $|\langle e, xy\rangle|=2$. Since $(xy)\varphi^{-1}\in E_S$, it follows from the preceding paragraph that either $(xy)\varphi^{-1}=ab^3a^2$ if $\mu=3$, or $(xy)\varphi^{-1}=b^2a^2$ if $\mu\in\{2,3\}$. In the former case, $|\langle e, xy\rangle|=|\langle ab, ab^3a^2\rangle|=|\{ab, ab^3a^2, ab^3a^3b\}|=3$; a contradiction. Now assume that $(xy)\varphi^{-1}=b^2a^2$. If $\mu=2$, then $|\langle e, xy\rangle|=|\langle ab, b^2a^2\rangle|=|\{ab, b^2a^2, b^2a^3b\}|=3$, and if $\mu=3$, then  $|\langle e, xy\rangle|=|\langle ab, b^2a^2\rangle|=|\{ab, b^2a^2, ab^3a^2,  b^2a^3b,ab^3a^3b \}|=5$, so we have a contradiction in each of these two situations. Therefore $xy\notin E_T$.

Since  $xy\notin E_T$ and $e(xy)=e=(xy)e$, we observe once again that $o(xy)=\infty$ and $e$ is an adjoined identity in the subsemigroup $\langle e, xy\rangle$ of $T$, so that, by \cite[Lemma 33.4]{key30}, $ab$ is either an adjoined identity or an adjoined zero in the subsemigroup $\langle ab, (xy)\varphi^{-1}\rangle$ of $S$. However, one can easily check that $ab$ is neither an identity nor a zero for any element of $L_{a^2}\setminus\{b^2a^2\}$ (when $\mu=2$) and of $L_{a^2}\setminus\{b^2a^2, ab^3a^2\}$ (when $\mu=3$). 

Since the assumption that $(xy)\varphi^{-1}\in L_{a^2}$ leads to a contradiction in all cases, it follows that  $(xy)\varphi^{-1}\notin L_{a^2}$ when $\nu=3$, and therefore, by symmetry, $(xy)\varphi^{-1}\notin R_{b^2}$ when $\mu=3$.
\vspace{0.05in}\\
\indent We have shown that none of the possibilities for $(xy)\varphi^{-1}$ listed in (\ref{eq:7th}) can actually occur. Since (\ref{eq:7th}) was deduced from the assumption that $e\neq xy$, we conclude that $e=xy$.  Finally, note that $(yx)^2=y(xy)x=yex=vu\in E_T\cap\cB(u,v)$. Hence $o(yx)\neq\infty$, from which it follows that $yx\in E_T$ and so $yx=(yx)^2=vu\neq uv=e$. This completes the proof of Lemma \ref{308}.     \epr
\setcounter{theorem}{12}

\bl\label{313} Suppose statement {\rm(II)} of Lemma {\rm\ref{303}} holds and $e=uv$. Then for all $k,l\in\N$ the following assertions are true: {\rm(i)} $x^ky^k=e$; {\rm(ii)} if $y^kx^l=y^mx^n$ for some $m,n\in\N$, then $k=m$ and $l=n$; {\rm(iii)} $y^kx^l=v^ku^l$ if statement {\rm(II)} of Lemma {\rm\ref{304}} holds, $y^kx^l=y^ku^l$ if $\varphi|_{L_{ab}}$ is an isomorphism of $L_{ab}$ onto $L_{ab}\Phi$, and in the latter case, if $\mu>1$, then $y^kx^l\neq v^ku^l$ if $k<\mu$; {\rm(iv)} $y^kx^k\in E_T\setminus\{e\}$ and $y^kx^l\notin E_T$ if $k\neq l$; and {\rm(v)} $y^kx^k>y^lx^l$ if $l>k$. 
\el
{\bf Proof.} (i) By Lemma \ref{308}, $x^ky^k=e$ when $k=1$. Suppose that $k\geq 2$. If $\varphi|_{L_{ab}}$ is an isomorphism of $L_{ab}$ onto $L_{ab}\Phi$, or if statement {\rm(II)} of Lemma {\rm\ref{304}} holds, then $v=ey$, so that $x^{k-1}ey^{k-1}=u^{k-1}v^{k-1}=e$. Therefore $x^ky^k=x^{k-1}(xy)y^{k-1}=x^{k-1}ey^{k-1}=e$. Assertion (i) will be used below without reference.
\vspace{0.06in}\\
\indent (ii) Suppose that $y^kx^l=y^mx^n$ for some $m,n\in\N$. Without loss of generality, we may assume that $k\leq m$. Then
 
\[e=(x^ky^k)(x^ly^l)=(x^ky^k)y^{m-k}(x^ny^l)=(ey^{m-k})(x^ny^l)=v^{m-k}(x^ny^l).\]
\vspace{0.02in}\\ 
If $n<l$, then $x^ny^l=(x^ny^n)y^{l-n}=v^{l-n}$, so that $e=v^{m-k}v^{l-n}=v^{(m-k)+(l-n)}$; a contradiction since $(m-k)+(l-n)>0$. Therefore $l\leq n$ and $x^ny^l=x^{n-l}(x^ly^l)=x^{n-l}e=u^{n-l}$. Hence $e=v^{m-k}u^{n-l}$ where $m-k\geq 0$ and $n-l\geq 0$, from which it follows that $k=m$ and $l=n$.
\vspace{0.05in}\\
\indent (iii) Suppose that statement {\rm(II)} of Lemma {\rm\ref{304}} holds. As shown in the proof of Lemma \ref{308}, $yx=vu$. Therefore 
\[y^kx^l=y^{k-1}(yx)x^{l-1}=y^{k-1}(vu)x^{l-1}=(y^{k-1}v)(ux^{l-1})=v^ku^l.\]

Now assume that $\varphi|_{L_{ab}}$ is an isomorphism of $L_{ab}$ onto $L_{ab}\Phi$. Since $y=ye$, it is immediate that $y^kx^l=(y^ke)x^l=y^k(ex^l)=y^ku^l$. Suppose that $\mu>1$ and $k<\mu$. If $y^kx^l=v^ku^l$, then $y^k=(y^kx^l)y^l=(v^ku^l)y^l=(v^kex^l)y^l=v^ke=v^k$; a contradiction. Hence $y^kx^l\neq v^ku^l$.
\vspace{0.05in}\\
\indent (iv) As shown in (iii), if statement {\rm(II)} of Lemma {\rm\ref{304}} is true, then $y^kx^k=v^ku^k$, so that $y^kx^k\in E_T\setminus\{e\}$. Suppose that $\varphi|_{L_{ab}}$ is an isomorphism of $L_{ab}$ onto $L_{ab}\Phi$. Since $ye=y$, we have $(y^kx^k)^2=y^k(x^ky^k)x^k=y^kex^k=y^kx^k$, so $y^kx^k\in E_T$. Note that $(yx)e=y(xe)=yu=yx$. In the proof of Lemma \ref{308} it was shown that $yx\in E_T\setminus\{e\}$. Assume that $k\geq 2$ and $y^kx^k=e$. Then $yx=(yx)(y^kx^k)=(yxy)y^{k-1}x^k=y^kx^k$, which contradict (ii). Therefore $y^kx^k\neq e$.

Suppose that $k\neq l$ and $y^kx^l\in E_T$. If $k<l$, then 
\[e=(x^ky^k)(x^ly^l)=x^k[(y^kx^l)(y^kx^l)]y^l=(x^ky^k)x^{l-k}(x^ky^k)(x^ly^l)=ex^{l-k}e=u^{l-k},\]
and if $k>l$, then
\[e=(x^ky^k)(x^ly^l)=x^k[(y^kx^l)(y^kx^l)]y^l=(x^ky^k)(x^ly^l)y^{k-l}(x^ly^l)=ey^{k-l}e=v^{k-l},\]
so we have a contradiction in both cases. It follows that $y^kx^l\notin E_T$ if $k\neq l$.
\vspace{0.05in}\\
\indent (v) Suppose that $l>k$. If statement (II) of Lemma \ref{304} holds, then since $\langle u,v\rangle=\cB(u,v)$, we conclude that $y^kx^k=v^ku^k>v^lu^l=y^lx^l$; note also that in this case $yx=vu<uv=e$.

Now assume that $\varphi|_{L_{ab}}$ is an isomorphism of $L_{ab}$ onto $L_{ab}\Phi$. Then
\[(y^kx^k)(y^lx^l)=y^k(x^ky^k)y^{l-k}x^l=y^key^{l-k}x^l=y^ky^{l-k}x^l=y^lx^l\]
and 
\[(y^lx^l)(y^kx^k)=y^lx^{l-k}(x^ky^k)x^k=y^lx^{l-k}ex^k=y^lx^{l-k}u^k=y^lu^l=y^lx^l,\]
that is, $y^kx^k\geq y^lx^l$. As shown in (ii), $y^kx^k\neq y^lx^l$. Hence $y^kx^k>y^lx^l$. Recall that if $\mu>1$, then $y\neq v$, and since $yx=vx$ implies $y=yxy=vxy=v$, it follows that $yx\neq vx=e(yx)$, so in this case $yx\nless e$. The proof of Lemma \ref{313} is complete.   \epr 
 
\bl\label{314} Suppose that statement {\rm(II)} of Lemma {\rm\ref{303}} holds and $e=uv$. Then 
\[T=\cB(u,v)\cup\langle x\rangle\cup\langle y\rangle\cup\{y^kx^l:k,l\in\N\}\;\text{ and }\;E_T= E_{\cB(u,v)}\cup\{y^kx^k: k\in\N\}.\]
\el
{\bf Proof.} Since $T=\langle x,y\rangle$, each element of $T$ can be written in the form $(x^{m_1}y^{n_1})\cdots(x^{m_r}y^{n_r})$ for some $r\in\N$ and some nonnegative integers $m_i$ and $n_i$ such that $m_i+n_i\geq 1$ for all $1\leq i\leq r$; as in \cite{key10}, we will call $(x^{m_1}y^{n_1}),\ldots,(x^{m_r}y^{n_r})$ the {\em syllables} of $(x^{m_1}y^{n_1})\cdots(x^{m_r}y^{n_r})$ and refer to $(x^{m_1}y^{n_1})\cdots(x^{m_r}y^{n_r})$ as an {\em $r$-syllable} $(x,y)$-word. Consider an arbitrary $1$-syllable word $x^my^n$. Using Lemma \ref{313}(i), we observe that
\setcounter{equation}{14}
\be
x^my^n=\left\{\!\! 
\begin{array}{l}x^m\\  y^n\\ x^my^m=e\\ x^{m-n}(x^ny^n)=x^{m-n}e=u^{m-n} \\ (x^my^m)y^{n-m}=ey^{n-m}=v^{n-m} \\ \end{array}
\begin{array}{l}\!\!\!\text{ if }n=0,\\ \!\!\!\text{ if }m=0,\\ \!\!\!\text{ if }m=n\geq 1, \\\!\!\!\text{ if }m>n\geq 1, \\ \!\!\!\text{ if }n>m\geq 1, 
				\end{array}\right.
\label{eq:8th}
\ee 
so $x^my^n\in\cB(u,v)\cup\langle x\rangle\cup\langle y\rangle$. Now take any $r>1$ such that the value of every $(r-1)$-syllable $(x,y)$-word is an element of $\cB(u,v)\cup\langle x\rangle\cup\langle y\rangle\cup\{y^kx^l:k,l\in\N\}$. Let $t$ be an arbitrary element of $T$ represented by an $r$-syllable $(x,y)$-word, say, $t=(x^{m_1}y^{n_1})\cdots(x^{m_{r-1}}y^{n_{r-1}})(x^my^n)$. Denote by $s$ the value of $(x^{m_1}y^{n_1})\cdots(x^{m_{r-1}}y^{n_{r-1}})$ in $T$. Thus $t=s(x^my^n)$ and, by our assumption, $s\in\cB(u,v)\cup\langle x\rangle\cup\langle y\rangle\cup\{y^kx^l:k,l\in\N\}$. If we show that for all possible values of $s$ we have $t\in\cB(u,v)\cup\langle x\rangle\cup\langle y\rangle\cup\{y^kx^l:k,l\in\N\}$, the equality $T=\cB(u,v)\cup\langle x\rangle\cup\langle y\rangle\cup\{y^kx^l:k,l\in\N\}$ will be established.
\vspace{0.06in}\\
\indent {\bf Case 1:} $s\in\cB(u,v)$.
\vspace{0.03in}\\ 
\indent Here $s=v^pu^q$ for some $p,q\geq 0$ such that $p+q\geq 1$. According to (\ref{eq:8th}), if $m, n\geq 1$, then $x^my^n\in\cB(u,v)$, in which case it is clear that $t=s(x^my^n)\in\cB(u,v)$. If $n=0$, then $x^my^n=x^m$, so that $t=v^pu^qx^m=v^pu^{q+m}\in\cB(u,v)$, and if $m=0$, then $x^my^n=y^n$ and hence $t=v^pu^qy^n=v^pu^qey^n=v^pu^qv^n\in\cB(u,v)$.
\vspace{0.06in}\\
\indent {\bf Case 2:} $s\in\langle x\rangle$.
\vspace{0.03in}\\
\indent Since $s=x^p$ for some $p\in\N$, we have $t=s(x^my^n)=x^p(x^my^n)=x^{p+m}y^n$ where $p+m>0$. Therefore, according to (\ref{eq:8th}), $t\in\cB(u,v)\cup\langle x\rangle$.
\vspace{0.06in}\\
\indent {\bf Case 3:} $s\in\langle y\rangle$.
\vspace{0.03in}\\
\indent In this case, $s=y^q$ for some $q\in\N$, so that $t=y^qx^my^n$. If $m=0$, then $t=y^{q+n}\in\langle y\rangle$, and if $n=0$, then $t=y^qx^m\in\{y^kx^l : k,l\in\N\}$. Suppose that $m,n\in\N$. If $m=n$, then $t=y^qe$, which shows that $t=y^q\in\langle y\rangle$ if $\varphi|_{L_{ab}}$ is an isomorphism of $L_{ab}$ onto $L_{ab}\Phi$, and $t=v^q\in\cB(u,v)$ if statement (II) of Lemma \ref{304} holds. By (\ref{eq:8th}), if $m>n$, then $t=y^qex^{m-n}$, so that $t=y^qx^{m-n}\in\{y^kx^l : k,l\in\N\}$ if $\varphi|_{L_{ab}}$ is an isomorphism of $L_{ab}$ onto $L_{ab}\Phi$, and $t=v^qu^{m-n}\in\cB(u,v)$ if statement (II) of Lemma \ref{304} is true. By (\ref{eq:8th}) again, if $n>m$, then $t=y^qv^{n-m}$ and therefore $t=y^{q+n-m}\in\langle y\rangle$ if $\varphi|_{L_{ab}}$ is an isomorphism of $L_{ab}$ onto $L_{ab}\Phi$, and $t=v^{q+n-m}\in\cB(u,v)$ if statement (II) of Lemma \ref{304} holds. The observations of this paragraph show that $t\in\cB(u,v)\cup\langle y\rangle\cup\{y^kx^l : k,l\in\N\}$.
\vspace{0.06in}\\
\indent {\bf Case 4:} $s\in\{y^kx^l:k,l\in\N\}$.
\vspace{0.03in}\\   
\indent Here $s=y^qx^p$ for some $p,q\in\N$ and hence $t=(y^qx^p)(x^my^n)=y^qx^{p+m}y^n$. Let $m'=p+m$. Then $m'\geq 1$ and $t=y^qx^{m'}y^n$, so using $m'$ instead of $m$ in the argument applied in Case 3 to $y^qx^my^n$ with $m\geq 1$, we conclude that $t\in\cB(u,v)\cup\langle y\rangle\cup\{y^kx^l : k,l\in\N\}$ in this case as well.    	
\vspace{0.06in}\\
\indent We have shown that $T=\cB(u,v)\cup\langle x\rangle\cup\langle y\rangle\cup\{y^kx^l:k,l\in\N\}$. By Lemma \ref{313}, if $k,l\in\N$, then $y^kx^l\in E_T$ if and only if $k=l$. Since no idempotent of $T$ lies in $\langle x\rangle\cup\langle y\rangle$, it follows that $E_T= E_{\cB(u,v)}\cup\{y^kx^k: k\in\N\}$. This completes the proof of Lemma \ref{314}.	\epr
\setcounter{theorem}{15}
\bl\label{316} In Lemmas {\rm\ref{303}} and {\rm\ref{304}} only statement {\rm(I)} is actually true.
\el
{\bf Proof.} As noted in the paragraph preceding Lemma \ref{308}, in view of duality, it will suffice to prove that statement (II) of Lemma \ref{303} does not hold under the assumption that $e=uv$.
\vspace{0.04in}\\
\indent Suppose, therefore, that $e=uv$ and statement (II) of Lemma \ref{303} is true. By Lemma \ref{314}, $T=\cB(u,v)\cup\langle x\rangle\cup\langle y\rangle\cup\{y^kx^l:k,l\in\N\}$ and $E_T= E_{\cB(u,v)}\cup\{y^kx^k: k\in\N\}$. Recall that $\nu\in\{2,3\}$. Then, according to Lemma \ref{210}, $\cB(a^2b,ab^2)\cap L_a=\emptyset$, and since $ba\in L_a\cap E_S$, it follows that $(ba)\varphi\in E_T\setminus E_{\cB(u,v)}$. 
\vspace{0.04in}\\
\indent If $\mu=1$, then $y=v$ and $\{y^kx^l : k,l\in\N\}=\{v^ku^l : k,l\in\N\}\subseteq\cB(u,v)$, which implies that $T=\cB(u,v)\cup\langle x\rangle\cup\langle y\rangle$ and $E_T=E_{\cB(u,v)}$, so that $(ba)\varphi\in E_T\setminus E_{\cB(u,v)}=\emptyset$; a contradiction. Hence $\mu>1$. Then $ab^2a\neq ba$ and $L_a\cap E_S=\{ab^2a, ba\}$; moreover,  $\cB(a^2b,ab^2)\cap R_b=\emptyset$ by Lemma \ref{210}. Since $\nu>1$, we also have $ba^2b\neq ba$ and $R_b\cap E_S=\{ba^2b, ba\}$. By Lemma \ref{207}(v), $\langle ab^2a, ba^2b\rangle=\{ab^2a^2b, ab^2a, ba, ba^2b\}$ is a four-element nonsingular rectangular band. Since $\cB(a^2b,ab^2)\cap L_a=\emptyset$ and $\cB(a^2b,ab^2)\cap R_b=\emptyset$, it follows that $(ab^2a)\varphi, (ba^2b)\varphi\in E_T\setminus E_{\cB(u,v)}$. By Lemma \ref{313}(iv), $(ab^2a)\varphi=y^kx^k$ and $(ba^2b)\varphi=y^lx^l$ for some distinct $k,l\in\N$, and therefore, in view of Lemma \ref{313}(v), $\langle y^kx^k, y^lx^l\rangle$ is a two-element semilattice. Thus we have a contradiction: $|\langle (ab^2a)\varphi, (ba^2b)\varphi\rangle|=2\neq 4=|\langle ab^2a, ba^2b\rangle|$. The proof of Lemma \ref{316} is complete.	\epr

\bl\label{317} Either both mappings $\varphi|_{R_{ab}}$ and $\varphi|_{L_{ab}}$ are isomorphisms or both are antiisomorphisms. If $\varphi|_{R_{ab}}$ and $\varphi|_{L_{ab}}$ are isomorphisms {\rm[}antiisomorphisms\,{\rm]}, then $\varphi|_{\cB(a^2b,\,ab^2)}$ is also an isomorphism {\rm[}antiisomorphism\,{\rm]} and $T=\langle x,y\rangle$ is a monogenic orthodox semigroup such that $e=xy=x^2y^2$ and $yx\neq y^2x^2$ {\rm[}$e=yx=y^2x^2$ and $xy\neq x^2y^2$\,{\rm]}; in fact, more precisely, $T=\cO_{(\nu,\,\mu)}(x,y)$ {\rm[}$T=\cO_{(\mu,\,\nu)}(y,x)${\rm]}. 
\el
{\bf Proof.} By Lemmas \ref{303}, \ref{304}, and \ref{316}, $\varphi|_{R_{ab}}$ is an isomorphism or an antiisomorphism of $R_{ab}$ onto $R_{ab}\Phi$, and $\varphi|_{L_{ab}}$ is an isomorphism or an antiisomorphism of $L_{ab}$ onto $L_{ab}\Phi$. Therefore, by Lemma \ref{307}, either both $\varphi|_{R_{ab}}$ and $\varphi|_{L_{ab}}$ are isomorphisms or they both are antiisomorphisms. 
\vspace{0.03in}\\
\indent Suppose that $\varphi|_{R_{ab}}$ and $\varphi|_{L_{ab}}$ are isomorphisms. Then, according to Lemma \ref{305}, $x=ex$, $u=xe$, $y=ye$, and $v=ey$, and hence, by Lemma \ref{306}, $e=uv=xey$ which implies that $\langle u,v\rangle=\cB(u,v)$ and $\varphi|_{\cB(a^2b,\,ab^2)}$ is an isomorphism of $\cB(a^2b,\,ab^2)$ onto $\cB(u,v)$. Applying Lemma \ref{302}, we conclude that $e=xy$ and $T=\langle x, y\rangle$ is a monogenic orthodox semigroup such that $xy=x^2y^2$ and $yx\neq y^2x^2$. 
\vspace{0.03in}\\
\indent By Lemma \ref{204}, $R_{ab}=\langle a,ab\rangle$, $L_{ab}=\langle b,ab\rangle$, $R_{xy}=\langle x,xy\rangle$, and $L_{xy}=\langle y,xy\rangle$. Since $\varphi|_{R_{ab}}$ is an isomorphism of $\langle a,ab\rangle$ onto $\langle x,xy\rangle$, it is clear that $a^k(ab)\neq a^k$ if and only if $x^k(xy)\neq x^k$ for $k\in\N$. Likewise, since $\varphi|_{L_{ab}}$ is an isomorphism of $\langle b,ab\rangle$ onto $\langle y,xy\rangle$, we have $(ab)b^l\neq b^l$ if and only if $(xy)y^l\neq y^l$ for $l\in\N$. It follows that $T=\cO_{(\nu,\,\mu)}(x,y)$.
\vspace{0.03in}\\
\indent The alternative statement holds by duality, so the proof of Lemma \ref{317} is complete.	\epr

\bl\label{318} If $\varphi|_{R_{ab}}$ and $\varphi|_{L_{ab}}$ are isomorphisms {\rm[}antiisomorphisms\,{\rm]} and if $\varphi(ba)=yx$ {\rm[}$\varphi(ba)=xy${\rm]}, then $(S\setminus R_b)\Phi = T\setminus R_y$ {\rm[}$(S\setminus R_b)\Phi = T\setminus L_y${\rm]} if $\mu>1$, and $(S\setminus L_a)\Phi = T\setminus L_x$ {\rm[}$(S\setminus L_a)\Phi = T\setminus R_x${\rm]} if $\nu>1$.  
\el
{\bf Proof.} Suppose that $\varphi|_{R_{ab}}$ and $\varphi|_{L_{ab}}$ are isomorphisms and $\varphi(ba)=yx$. By Lemma \ref{317}, $T=\cO_{(\nu,\,\mu)}(x,y)$. Assume that $\mu>1$. By Lemma \ref{205}(i), $S\setminus R_b\leq S$ and $T\setminus R_y\leq T$. Suppose that $(S\setminus R_b)\Phi\cap R_y\neq\emptyset$. If $y\in(S\setminus R_b)\Phi$, then using the fact that $x=a\varphi\in(S\setminus R_b)\Phi$, we conclude that $T=\langle x,y\rangle\subseteq(S\setminus R_b)\Phi\subset T$; a contradiction. Thus $y\notin(S\setminus R_b)\Phi$. Since $ba\in R_b$ and $(ba)\varphi=yx$, it is immediate that $yx\notin(S\setminus R_b)\Phi$. If $yx^k\in(S\setminus R_b)\Phi$ for some $k\geq 2$, then $y^k=(b^k)\varphi\in(S\setminus R_b)\Phi$ whence $y=y(xy)=y(x^ky^k)=(yx^k)y^k\in(S\setminus R_b)\Phi$; a contradiction. Therefore $yx^k\notin(S\setminus R_b)\Phi$ for all $k\in\N$. Finally, if $yx^ky\in(S\setminus R_b)\Phi$ for some $k\geq 2$, then $yx^k=(yx^ky)x\in(S\setminus R_b)\Phi$, contrary to the preceding conclusion. Thus $yx^ky\notin(S\setminus R_b)\Phi$ for all $k\geq 2$. We have shown that $(S\setminus R_b)\Phi\cap R_y=\emptyset$, that is, $(S\setminus R_b)\Phi\subseteq T\setminus R_y$. Using $\varphi^{-1}$ and $\Phi^{-1}$ instead of $\varphi$ and $\Phi$, respectively, we obtain $(T\setminus R_y)\Phi^{-1}\subseteq S\setminus R_b$. Therefore $(S\setminus R_b)\Phi = T\setminus R_y$. If $\nu>1$, then $(S\setminus L_a)\Phi = T\setminus L_x$ is deduced by a symmetric argument. Since the alternative statement holds by duality, the proof of Lemma \ref{318} is complete. \epr

\bl\label{319} If $\varphi|_{R_{ab}}$ and $\varphi|_{L_{ab}}$ are isomorphisms {\rm[}antiisomorphisms\,{\rm]}, then $\varphi$ is an isomorphism {\rm[}antiisomorphism\,{\rm]} of $S$ onto $T$.
\el
{\bf Proof.} Suppose that $\varphi|_{R_{ab}}$ and $\varphi|_{L_{ab}}$ are isomorphisms. Then, by Lemma \ref{317}, $\varphi|_{\cB(a^2b,\,ab^2)}$ is an isomorphism of $\cB(a^2b,\,ab^2)$ onto $\cB(x^2y,xy^2)$ and $T=\cO_{(\nu,\,\mu)}(x,y)$. Recall that every element of $S$ [of $T$] is represented by a unique reduced $(a,b)$-word $a^ib^ma^nb^j$ [$(x,y)$-word $x^iy^mx^ny^j$] where $m-1<\mu$ and $n-1<\nu$; {\em in what follows, whenever we say that $a^ib^ma^nb^j$ {\rm[}$x^iy^mx^ny^j${\rm]} is a reduced $(a,b)$-word {\rm[}$(x,y)$-word\,{\rm]}, these restrictions on $m$ and $n$ will be assumed and often used without mention}. Note that an abridged $(a,b)$-word $a^ib^ma^nb^j$ is reduced precisely when the corresponding $(x,y)$-word $x^iy^mx^ny^j$ is reduced. It is also clear that $\varphi$ is an isomorphism of $S$ onto $T$ if and only if $(a^ib^ma^nb^j)\varphi=x^iy^mx^ny^j$ for each reduced $(a,b)$-word $a^ib^ma^nb^j$. 
\vspace{0.05in}\\
\indent {\bf Case 1:} $\mu>1$ and $\nu>1$.
\vspace{0.03in}\\ 
\indent In this case, by Lemma \ref{208}(i),
\setcounter{equation}{19}
\be
S=\cB(a^2b,ab^2)\cup\cB(ba^3b,b^2a^2b)\cup\cB(ba^2,b^2a)\cup\cB(ab^2a^2,ab^3a)\cup\langle a\rangle\cup\langle b\rangle
\label{eq:9th}
\ee
and 
\be
T=\cB(x^2y,xy^2)\cup\cB(yx^3y,y^2x^2y)\cup\cB(yx^2,y^2x)\cup\cB(xy^2x^2,xy^3x)\cup\langle x\rangle\cup\langle y\rangle.
\label{eq:10th}
\ee
Since $\varphi|_{R_{ab}}$ and $\varphi|_{L_{ab}}$ are isomorphisms, we have $xy^2\cL xy\cR x$ whence $xy^2\cR xy^2x\cL x$ by \cite[Theorem 2.17]{key5}. From the fact that $\varphi|_{\cB(a^2b,\,ab^2)}$ is an isomorphism of $\cB(a^2b,\,ab^2)$ onto $\cB(x^2y,xy^2)$, it follows that $(a^ib^ma^nb^j)\varphi=x^iy^mx^ny^j$ for each reduced $(a,b)$-word $a^ib^ma^nb^j\in\cB(a^2b,\,ab^2)$. In particular, $(ab^2a^2b)\varphi=xy^2x^2y$ and $xy^2x^2y\cR xy^2$. Therefore $xy^2x^2y\cR xy^2x$. By Lemma \ref{207}(v), $\langle ab^2a^2b, ba\rangle=\{ab^2a^2b, ab^2a, ba, ba^2b\}$ and $\langle xy^2x^2y, yx\rangle=\{xy^2x^2y, xy^2x, yx, yx^2y\}$ are four-element nonsingular rectangular bands (which coincide, respectively, with $D^{E_S}_{ab^2a^2b}$ and $D^{E_T}_{xy^2x^2y}$). According to Result \ref{105}, $\varphi$ is an isomorphism or an antiisomorphism of $\langle ab^2a^2b, ba\rangle$ onto $\langle (ab^2a^2b)\varphi, (ba)\varphi\rangle$. Since $(ab^2a^2b)\varphi=xy^2x^2y$ and $D^{E_T}_{xy^2x^2y}$ is the only nonsingular rectangular band containing $xy^2x^2y$, it follows that $\langle (ab^2a^2b)\varphi, (ba)\varphi\rangle=\langle xy^2x^2y, yx\rangle$ and $(ba)\varphi=yx$, so $\{ab^2a,ba^2b\}\varphi=\{xy^2x,yx^2y\}$.
\vspace{0.04in}\\
\indent By Lemma \ref{318}, $(S\setminus R_b)\Phi = T\setminus R_y$ since $\mu>1$ and $(ba)\varphi=yx$. Clearly, $ab^2a\in S\setminus R_b$, $xy^2x\in T\setminus R_y$, and $yx^2y\in R_y$, so using the fact that $\{ab^2a,ba^2b\}\varphi=\{xy^2x,yx^2y\}$, we conclude that $(ab^2a)\varphi=xy^2x$ and $(ba^2b)\varphi=yx^2y$. Thus, in view of Lemma \ref{207}${\rm(ii, iii, iv)}$, $\varphi$ maps the identity elements of $\cB(ab^2a^2,ab^3a)$, $\cB(ba^2,b^2a)$, and $\cB(ba^3b,b^2a^2b)$ onto the identities of $\cB(xy^2x^2,xy^3x)$, $\cB(yx^2,y^2x)$, and $\cB(yx^3y,y^2x^2y)$, respectively. By Result \ref{104} and Lemma \ref{208}(ii), we have $\cB(ab^2a^2,ab^3a)\Phi\subseteq\cB(xy^2x^2,xy^3x)$, $\cB(ba^2,b^2a)\Phi\subseteq\cB(yx^2,y^2x)$, and $\cB(ba^3b,b^2a^2b)\Phi\subseteq\cB(yx^3y,y^2x^2y)$. By a symmetric argument, using $\varphi^{-1}$ and $\Phi^{-1}$ instead of $\varphi$ and $\Phi$, respectively, we obtain the converse inclusions. It follows that $\cB(ab^2a^2,ab^3a)\Phi=\cB(xy^2x^2,xy^3x)$, $\cB(ba^2,b^2a)\Phi=\cB(yx^2,y^2x)$, and $\cB(ba^3b,b^2a^2b)\Phi=\cB(yx^3y,y^2x^2y)$. Moreover, by Result \ref{104}, $\varphi|_{\cB(ab^2a^2,ab^3a)}$ [$\varphi|_{\cB(ba^2,b^2a)}$, $\varphi|_{\cB(ba^3b,b^2a^2b)}$] is an isomorphism or an antiisomorphism of $\cB(ab^2a^2,ab^3a)$ onto $\cB(xy^2x^2,xy^3x)$ [of $\cB(ba^2,b^2a)$ onto $\cB(yx^2,y^2x)$, of $\cB(ba^3b,b^2a^2b)$ onto $\cB(yx^3y,y^2x^2y)$]. 
\vspace{0.04in}\\
\indent Suppose that $\varphi|_{\cB(ab^2a^2,ab^3a)}$ is an antiisomorphism of $\cB(ab^2a^2,ab^3a)$ onto $\cB(xy^2x^2,xy^3x)$. Then $(ab^2a^2)\varphi=xy^3x$ and $(ab^3a)\varphi=xy^2x^2$. By Lemma \ref{318}, $(S\setminus L_a)\Phi = T\setminus L_x$ since $\nu>1$. Clearly, $ab^2a^2\in S\setminus L_a$ and $xy^3x\in L_x$. Hence $(ab^2a^2)\varphi\in(S\setminus L_a)\Phi\cap L_x=\emptyset$; a contradiction. Therefore $\varphi|_{\cB(ab^2a^2,ab^3a)}$ is an isomorphism of $\cB(ab^2a^2,ab^3a)$ onto $\cB(xy^2x^2,xy^3x)$, so that $(a^ib^ma^nb^j)\varphi=x^iy^mx^ny^j$ for each reduced $(a,b)$-word $a^ib^ma^nb^j\in\cB(ab^2a^2,ab^3a)$. Similarly, one can show that $\varphi|_{\cB(ba^2,\, b^2a)}$ is an isomorphism of $\cB(ba^2,b^2a)$ onto $\cB(yx^2,y^2x)$, and $\varphi|_{\cB(ba^3b,\, b^2a^2b)}$ is an isomorphism of $\cB(ba^3b,b^2a^2b)$ onto $\cB(yx^3y,y^2x^2y)$, and hence $(a^ib^ma^nb^j)\varphi=x^iy^mx^ny^j$ for each reduced $(a,b)$-word $a^ib^ma^nb^j$ in $\cB(ba^2,b^2a)$ or in $\cB(ba^3b,b^2a^2b)$, respectively. Since $\varphi|_{\langle a\rangle}$ is an isomorphism of $\langle a\rangle$ onto $\langle x\rangle$, and $\varphi|_{\langle b\rangle}$ is an isomorphism of $\langle b\rangle$ onto $\langle y\rangle$, in view of the above and formulas (\ref{eq:9th}) and (\ref{eq:10th}), it follows that $\varphi$ is an isomorphism of $S$ onto $T$.
\vspace{0.05in}\\
\indent {\bf Case 2:} $\mu=1$.
\vspace{0.03in}\\
\indent In this case, by Lemma \ref{209}(i),
\be 
S=\cB(a^2b,b)\cup\cB(ba^2,b^2a)\cup\langle a\rangle
\label{eq:11th}
\ee
and
\be 
T=\cB(x^2y,y)\cup\cB(yx^2,y^2x)\cup\langle x\rangle.
\label{eq:12th}
\ee
Since $\varphi|_{\langle a, ab\rangle}\colon \langle a, ab\rangle \to \langle x, xy\rangle$ and $\varphi|_{\cB(a^2b,b)}\colon \cB(a^2b, b) \to \cB(x^2y, y)$ are isomorphisms, we conclude that $(a^k)\varphi=x^k$ for all $k\in\N$, and $(a^nb)\varphi=x^ny$ and $(b^ma^nb)\varphi=y^mx^ny$ for all $m,n\in\N$ such that $n-1<\nu$, so $(ba^2b)\varphi=yx^2y$. Note that $\langle ba^2b, ba\rangle$ and $\langle yx^2y, yx\rangle$ are two-element right singular semigroups. Hence $ba\notin\cB(a^2b,b)$ and $yx\notin\cB(x^2y,y)$, which implies that $(ba)\varphi\notin\cB(x^2y, y)$. It is clear that if a semilattice $(\subset T)$ has $yx^2y$ as one of its elements, then that semilattice is contained in $\cB(x^2y, y)$. It follows that $\langle ba^2b, ba\rangle\Phi=\langle yx^2y,(ba)\varphi\rangle$ is a two-element singular semigroup, and the only such subsemigroup of $T$ is the two-element right singular semigroup $\langle yx^2y, yx\rangle$. Therefore $(ba)\varphi=yx$. Since $ba$ and $yx$ are identity elements of $\cB(ba^2,b^2a)$ and $\cB(yx^2,y^2x)$, respectively, in view of Result \ref{104} and Lemma \ref{209}(ii), $\cB(ba^2,b^2a)\Phi\subseteq\cB(yx^2,y^2x)$. By a symmetric argument, using $\varphi^{-1}$ and $\Phi^{-1}$ instead of $\varphi$ and $\Phi$, respectively, we have the converse inclusion. Thus $\cB(ba^2,b^2a)\Phi=\cB(yx^2,y^2x)$ and, by Result \ref{104}, $\varphi|_{\cB(ba^2,b^2a)}$ is an isomorphism or an antiisomorphism of $\cB(ba^2,b^2a)$ onto $\cB(yx^2,y^2x)$.
\vspace{0.04in}\\
\indent Suppose that $\varphi|_{\cB(ba^2,b^2a)}$ is an antiisomorphism of $\cB(ba^2,b^2a)$ onto $\cB(yx^2,y^2x)$, in which case $(ba^2)\varphi=y^2x$ and $(b^2a)\varphi=yx^2$. By Lemma \ref{318}, $(S\setminus L_a)\Phi = T\setminus L_x$ since $\nu>1$ and $(ba)\varphi=yx$. Note that $ba^2\in S\setminus L_a$ and $y^2x\in L_x$. Hence $(ba^2)\varphi\in(S\setminus L_a)\Phi\cap L_x=(T\setminus L_x)\cap L_x$; a contradiction. Therefore $\varphi|_{\cB(ba^2,b^2a)}$ is an isomorphism of $\cB(ba^2,b^2a)$ onto $\cB(yx^2,y^2x)$, so that $(b^ma^n)\varphi=y^mx^n$ for all $m,n\in\N$. In view of the above and formulas (\ref{eq:11th}) and (\ref{eq:12th}), it follows that $\varphi$ is an isomorphism of $S$ onto $T$.   
\vspace{0.05in}\\
\indent {\bf Case 3:} $\nu=1$.
\vspace{0.03in}\\
\indent In this case, $\varphi$ is an isomorphism of $S$ onto $T$ by symmetry with Case 2.
\vspace{0.08in}\\
\indent We have shown that if $\varphi|_{R_{ab}}$ and $\varphi|_{L_{ab}}$ are isomorphisms, then in all possible cases $\varphi$ is an isomorphism of $S$ onto $T$. To finish the proof of Lemma \ref{319}, it remains to note that if $\varphi|_{R_{ab}}$ and $\varphi|_{L_{ab}}$ are antiisomorphisms, then $\varphi$ is an antiisomorphism of $S$ onto $T$ by the dual of the above argument. \epr
\vspace{0.1in}\\
\indent From Lemmas \ref{317} and \ref{319}, it follows that $\varphi$ is an isomorphism or an antiisomorphism of $S$ onto $T$. Thus $S$ is strongly lattice determined. The proof of Theorem \ref{301} is complete. \epr
\vspace{0.15in}\\
\begin{center}
{\bf Acknowledgement}
\end{center}
\sm Part of the work on this paper was done while I had an appointment as a Visiting Professor at the Department of Mathematics at Stanford University from January 1 through June 30 of 2009. This support is greatly appreciated. I also would like to thank the department's faculty and staff for their hospitality and for creating excellent working conditions.
\vspace{0.15in}\\

\end{document}